\documentclass[main]{siamart220329}

\usepackage{amsfonts,subfigure}
\usepackage{graphicx}
\usepackage{epstopdf}
\usepackage{algorithmic}
\usepackage{hyperref}
 \hypersetup{
     colorlinks=true,
     linkcolor=black,
     filecolor=black,
     citecolor = black,      
     urlcolor=black,
     }
\usepackage{amsmath, amssymb}
\usepackage{relsize}
\usepackage{pifont}
\usepackage{enumitem}
\usepackage{oubraces}
\usepackage{listings}
\usepackage{booktabs}
\usepackage{xfrac}
\usepackage{cleveref}
\usepackage{caption}
\usepackage{tabularx}
\usepackage{array}
\usepackage{multirow}
\usepackage{mathrsfs}
\usepackage{microtype}
\usepackage[textfont=it,belowskip=-15pt]{caption}

\usepackage{tikz}
\usetikzlibrary{arrows,backgrounds,patterns}
\usetikzlibrary{matrix,shapes,fit,calc,shadows,plotmarks}
\usetikzlibrary{decorations}
\usetikzlibrary{decorations.markings}
\RequirePackage{pgfplotstable}
\RequirePackage{pgfplots}
\usepackage{pgfplotstable}
\usepackage{xparse}
\pgfplotsset{compat=1.18} 

\newcommand\eps{\epsilon}
\renewcommand*{\i}{\mathrm{i}}
\newcommand*{\e}{\mathrm{e}}
\newcommand*{\im}{\mathrm{i}}
\newcommand*{\ep}{\epsilon}
\newcommand*{\de}{\operatorname{d\!}{}} 
\newcommand{\dd}[2]{\frac{\de#1}{\de#2}}
\newcommand{\pd}[2]{\frac{\partial#1}{\partial#2}}

\newcommand{\myred}[1]{{\color{black}#1}}
\def\Xint#1{\mathchoice
   {\XXint\displaystyle\textstyle{#1}}%
   {\XXint\textstyle\scriptstyle{#1}}%
   {\XXint\scriptstyle\scriptscriptstyle{#1}}%
   {\XXint\scriptscriptstyle\scriptscriptstyle{#1}}%
   \!\int}
\def\XXint#1#2#3{{\setbox0=\hbox{$#1{#2#3}{\int}$}
     \vcenter{\hbox{$#2#3$}}\kern-.5\wd0}}

\def\dashint{\Xint-}

\headers{J. Shelton, S. Crew, P. H. Trinh}{The higher-order Stokes phenomenon}
\title{Exponential asymptotics and higher-order Stokes phenomenon in singularly perturbed ODEs
\thanks{Submitted to the editors 01/03/2023. JS and SC are joint lead authors.}
\funding{PHT is supported by the EPSRC grant EP/V012479/1. JS is supported by the EPSRC grant EP/W522491/1. This work was supported by EPSRC grant EP/R014604/1.}
}

\author{Josh Shelton\thanks{Department of Mathematical Sciences, University of Bath, BA2 7AY, UK (\email{jas45@st-andrews.ac.uk},\email{p.trinh@bath.ac.uk}).}
\and Samuel Crew\thanks{Department of Mathematics, Imperial College London, 180 Queen’s Gate, London SW7 2AZ, UK. (\email{samuel.c.crew@gmail.com})}
\and Philippe H. Trinh\footnotemark[2]}

\AtBeginDocument{\let\latexlabel\label}
\begin{document}

\maketitle
\begin{abstract} 
The higher-order Stokes phenomenon can emerge in the asymptotic analysis of many problems governed by singular perturbations.
Indeed, over the last two decades, the phenomena has appeared in many physical applications, from  acoustic and optical wave phenomena and gravity-capillary ripples, to models of crystal growth and equatorial Kelvin waves.
It emerges in a generic fashion in the exponential asymptotics of higher-order ordinary and partial differential equations.
The intention of this work is to highlight its importance, and develop further practical methodologies for the study of higher-order Stokes phenomena, primarily for general non-integrable problems.
Our formal methodology is demonstrated through application to a second-order linear inhomogeneous ODE that exemplifies the simplest example of higher-order Stokes phenomena.
In this model problem, the Borel transform can be derived explicitly, and this gives insight into the beyond-all-orders structure.
We review and study additional examples, with physically-important connections, including higher-order ODEs and eigenvalue problems.
\end{abstract}


\section{Introduction} \label{sec:intro}

Outside the exponential asymptotics community, the concept of the \emph{higher-order Stokes phenomenon} (HOSP) is still considered to be relatively obscure. Yet,  as we shall explain in this work, it occurs generically and routinely in divergent asymptotic expansions. Deep understanding of the phenomenon seems to be necessary in order to perform beyond-all-orders analysis of multi-dimensional problems (such as in partial differential equations) arising in many physical applications \cite{howls_2004,chapman_2005,body2005exponential}. The emphasis of this work is to provide some relatively simple examples of the HOSP in the context of singularly-perturbed ordinary differential equations. For example, we show that it occurs in the asymptotic analysis of the seemingly innocuous equation,
\begin{equation} \label{eq:themainone_intro}
 \eps^2 y''+ \eps (1+z)y' + z y = 1
\end{equation}
with $\ep \to 0$ and $y = y(z)$. We shall explain why the above is likely the simplest example of HOSP, and we develop applicable techniques for its study. Finally, our work serves to survey perspectives of the phenomenon across disciplines in asymptotic analysis.

We first review the regular Stokes phenomenon. In applications to singularly-perturbed ordinary differential equations, the solution is typically asymptotically expanded as
\begin{equation} \label{eq:base}
y(z; \ep) \sim  y_0(z) + \ep y_1(z) + \cdots + \ep^n y_n(z) + \cdots,
\end{equation}
where $\ep \to 0$ and $z \in \mathbb{C}$, and refer to such an expansion as the {\itshape base series}. The late terms, $y_n$, diverge as $n \to \infty$ and this divergence is coupled to the existence of Stokes lines in $z\in\mathbb{C}$. Across such lines, exponentially-small terms, say of $O(\e^{-\chi(z)/\ep})$, are switched-on within a small boundary layer via the (regular) Stokes phenomenon. By now, such ideas are well established (cf. reviews in e.g. Dingle \cite{dingle_book} and Berry \cite{berry_1989}), and there exists a plethora of tools in asymptotics beyond-all-orders to study these effects in the context of ordinary and partial differential equations, difference equations, integral equations, and beyond.

Typically, the divergence of $y_n$ is approximated by a factorial-over-power growth as $n \to \infty$. However, this is generally only the leading-order term in an asymptotic expansion in decreasing powers of $n$. If $1/n$ is now considered to be the perturbative parameter, then it is interesting to consider if $y_n$ can, itself, exhibit Stokes phenomena when analytically continued in $z$. Indeed this is the higher-order Stokes phenomenon associated with {\itshape higher-order Stokes lines} (HOSL). The standard Stokes phenomenon described previously is linked to the divergence of $y_n$, but now, since these approximations may be non-uniform in $\mathbb{C}_{z}$, the regular Stokes line structure is consequently altered. 

For the practitioner interested in the prediction of exponentially-small terms, two significant changes can occur as a consequence of HOSP:
\begin{enumerate}[label=(\roman*),leftmargin=*, align = left, labelsep=\parindent, topsep=3pt, itemsep=2pt,itemindent=0pt]
\item Stokes lines can be truncated or there can be new Stokes lines; and 
\item there can be atypical values of the Stokes multiplier(s). 
\end{enumerate}
We shall explore both of these consequences in this work.

\subsection{Prior literature and ongoing work}
There are a number of seminal works that study HOSP, though not necessarily in the same context of our work here. The origins of the phenomenon are generally attributed to the work of Berk \emph{et al.} \cite{berk1982new}. In their paper, titled ``\emph{New Stokes' line in WKB theory}", the altered structure of Stokes lines was discovered in the context of an integral equation (Pearcey's integral). The theory of such ``new Stokes lines" was then extensively studied in the exact WKBJ (Liouville-Green) communities (cf. Aoki~\cite{aoki1994new} and references therein). Further WKBJ examples with exact integral solutions were studied by Aoki \emph{et al.} \cite{aoki1998exact, aoki2001exact}; a singularly perturbed Painlev\'e equation appears in the work of Honda~\cite{honda2007stokes}; and HOSP appears as well in the analysis of the Henon map in quantum chaos by Shudo~\cite{shudo2007role}. We refer readers to reviews by \emph{e.g.} Aoki \emph{et al.} \cite{aoki2001exact}, Honda \emph{et al.} \cite{honda2015virtual}, and Takei \cite{takei2017wkb} and references therein for a historical survey. 

We note an important line of research in the early 2000s, initiated by a triplet of works by Howls \emph{et al.} \cite{howls_2004}, Chapman and Mortimer \cite{chapman_2005}, and Body \emph{et al.} \cite{body2005exponential}, which characterise some important connections between beyond-all-orders analysis of partial differential equations and the HOSP. These three works apply a number of different methodologies to beyond-all-order asymptotics, and present complementary views of similar problems. In particular, the work of Howls \emph{et al.} \cite{howls_2004} returned to analyse the Pearcey integral, discussed above, and established a connection between HOSP and the prior theory of hyperterminants developed by Olde Daalhuis \cite{daalhuis1996hyperterminants, daalhuis1998hyperterminants}. A particularly simple example of the HOSP was also presented by Olde Daalhuis \cite{daalhuis2004higher} in the context of an inhomogeneous linear second-order ODE. 

There have been a number of recent exciting works concerning HOSP in both applications and theoretical studies. For example, we highlight the work of Nemes \cite{nemes2022dingle} who studied atypical Stokes multipliers arising in the asymptotic representation of the Gamma function. Further applications have been discovered in fluid dynamics by Lustri \emph{et al.} \cite{lustri2019three} for three dimensional water waves, and Shelton {et al.} \cite{shelton2022Hermite,shelton2023kelvin} for geophysical instabilities.
 Most recently, we note unpublished investigations by Olde Daalhuis \cite{smoothingAdri} into the smoothing of Stokes multipliers of an $\ep$ expansion, and King \cite{iniKing} on the occurrence of HOSP in linear PDEs.

\subsection{Outline of our work}

\myred{We begin in \S\ref{sec:modelprob} by examining our minimal example \eqref{eq:themainone_intro} in the Borel plane $\mathbb{C}_{w}$.
An exact solution is found, which yields explicit insight into the connection between multiple singularities in the Borel plane and the higher-order Stokes phenomenon.
We then consider the same example in \S\ref{sec:physicalhospsec} with respect to asymptotics of the solution in the physical plane $\mathbb{C}_z$.} It is shown that a new divergent series (in the form of a $1/n$ expansion) can emerge within the late-terms themselves, which we denote the {\it late-late-terms} of the expansion. Borel summation applied to the tail of this divergent series reveals the HOSP, in which new components of the late-terms are seen to smoothly switch on with an error function dependence across higher-order Stokes lines.
Finally, in \S\ref{sec:otherexamples} we apply our general techniques to detect the higher-order Stokes lines in three previously encountered ODEs. These problems had previously been studied through methods that cannot be easily generalised to arbitrary nonlinear differential equations---such as steepest descent analysis on an integral representation, or evaluation of a known recurrence relation. 

\section{The model problem}\label{sec:modelprob}
We now introduce a new model example whose Stokes line structure contains the simplest possible realisation of the higher-order Stokes phenomenon, due to the presence of one branch point and a simple pole in the Borel plane.
This example is given by the following singularly-perturbed second-order linear inhomogeneous differential equation
\begin{subequations}
\begin{equation}\label{eq:themainone}
    \eps^2 y''+ \eps (1+z)y' + z y = 1,
    \end{equation}
where $y(z,\ep)$ is the unknown function and $0<\ep \ll 1$ is a small parameter. 
\myred{We will explore the Stokes line structure of this example throughout $z \in \mathbb{C}$ subject to the behavioural condition
\begin{equation}\label{eq:themainoneBC}
y(z) \sim \frac{1}{z} \quad \text{as} \quad \lvert z \rvert \to \infty,~ \text{arg}[z]=-2 \pi /3.
\end{equation}
Note that as $\lvert z \rvert \to \infty$ the far-field solution of \eqref{eq:themainone} is of the form $y \sim z^{-1} + \alpha z^{-1} \mathrm{e}^{-z^2/(2 \epsilon)}+\beta \mathrm{e}^{-z/\ep}$. Along the limit specified in condition \eqref{eq:themainoneBC} both exponentials in this far-field solution are exponentially dominant, hence requiring $\alpha=0$ and $\beta=0$.
Thus, the behavioural condition \eqref{eq:themainoneBC} provides two constraints for our second-order ODE.
In terms of the Stokes line structure and asymptotic transseries, this will correspond to only the base expansion \eqref{eq:base} existing in the lower-half plane, $\text{Im}[z]<0$.

}
\end{subequations}
\subsection{The Borel plane $\mathbb{C}_{w}$}\label{sec:mainborel}
Borel resummation encodes the parametric transseries (all algebraic and exponential orders of the solution expansion for $y(z,\epsilon)$ as $\epsilon \to 0$) as a single function, $y_B(w,z)$, which is holomorphic in $w$.
To see this, consider the formal base series \eqref{eq:base}, $y(z,\ep)\sim \sum_{n=0}^{\infty}\ep^n y_n(z)$. The Borel transform with respect to $\epsilon$ is a power series on the \textit{Borel plane}, $w \in \mathbb{C}_w$, defined by
\begin{equation}\label{eq:Borelexpansion}
    \myred{y_B(w,z) := \sum_{n=0}^{\infty}w^n \frac{y_{n}(z)}{\Gamma(n+1)}.}
\end{equation}
We assume throughout this work that $y_B$ is endlessly analytically continuable (cf. Mitschi \emph{et al.} \cite{mitschi2016divergent}) except at a finite number of points $z \in \mathbb{C}_z$, and extends to a function with at most exponential growth along any ray emanating from the origin. Under these conditions we may then define the inverse Borel transform of $y_B$ via 
\begin{equation}\label{eq:InverseBorel}
    y(z; \ep)=\mathcal{B}^{-1}[y_B]:=\frac{1}{\ep}\int_0^{\infty}y_B(w,z)\e^{-w/\ep}\de{w},
\end{equation}
which expresses the solution in the form of a Laplace transform. \myred{Note that our definitions above may differ from other treatments, e.g. Dorigoni \cite{dorigoni2019introduction}, by the $1/\ep$ scaling factor or shift of the indices. In this work, we will primarily consider $\ep$ real and positive; the analytic structure of $\ep \in \mathbb{C}$ may be important in certain problems as well, and consideration of this will require deformation of the contour in \eqref{eq:InverseBorel}.}

Analysis of the Borel transform, $y_B$, provides deep insight into the Stokes phenomenon, HOSP, and many of the calculations we will present later in this section. Remarkably, in the case of the model problem \eqref{eq:themainone}, the Borel transform can be developed explicitly. In order to solve for $y_B$, the integral representation \eqref{eq:InverseBorel} can be substituted into \eqref{eq:themainone} to develop the necessary equation for $y_B$ \cite{dorigoni2019introduction}. The Borel transform is governed by the PDE
\begin{equation}\label{eq:BorelPDE}
\pd{^2 y_B}{z^2} +(1+z)\myred{\pd{^2 y_B}{z \partial w}} + z \pd{^2 y_B}{w^2} = 0, 
\end{equation}
\myred{with boundary conditions $y_B(0, z) = y_0(z)$ and $\partial_{w}y_B(0,z)=y_1(z)$}. In using $y_0=1/z$ and $y_1=(1+z)/z^3$, equation \eqref{eq:BorelPDE} is found to have the exact solution
\begin{equation}\label{eq:exactBorel}
y_B(w,z) = \frac{1}{2(\chi_2-w)} \bigg(1+\frac{z-1}{\sqrt{2(\chi_1-w)}}\bigg),
\end{equation}
where we have denoted $\chi_1=z^2/2$ and $\chi_2=z-1/2$. It can subsequently be verified that expanding \eqref{eq:exactBorel} as $w \to 0$ recovers each order of the asymptotic expansion \eqref{eq:base} via definition \eqref{eq:Borelexpansion}, \myred{\emph{i.e.} $y_B = y_0(z)/\Gamma(1) + [y_1(z)/\Gamma(2)]w + [y_2(z)/\Gamma(3)]w^2 + \cdots$.} The key features of \eqref{eq:exactBorel} worth noting are the pole at $w=\chi_2$ and the branch point at $w=\chi_1$. It is both of these that will lead to the HOSP, and this is what motivated our minimal example \eqref{eq:themainone}. This seems to be the simplest arrangement of singularities that both displays HOSP and has an exact Borel solution. \myred{We note that, in the practical study of solutions to differential equations, particularly nonlinear cases, it is not typically possible to find such a closed form solution to the Borel PDE.}

\subsection{Stokes phenomenon in the Borel plane}
We begin by explaining how the regular Stokes phenomenon can be understood in the Borel plane.
Consider the asymptotic evaluation of the integral \eqref{eq:InverseBorel} for small real-valued $\epsilon$.
Provided singularities of $y_B$ do not lie on the integration contour, we may demonstrate using the integral formula of the Gamma function that $\mathcal{B}^{-1}[y_B] = y_0(z) + \epsilon y_1(z) + \cdots$. However, if there is an  algebraic singularity of $y_B$ at $w=\chi_1(z)$, where 
\begin{equation}
    y_B(w,z) \sim \frac{y_0^{(1)}(z)}{(\chi_1(z)-w)^{\alpha_1}} \quad \text{as}~ w \to \chi_{1}(z),
\end{equation}
and the singularity, under variation of $z$, crosses the integration contour $(0,\infty)$ of \eqref{eq:InverseBorel}, we pick up an extra contribution to the asymptotics of $\mathcal{B}^{-1}[y_B]$. This contribution is from the Hankel contour integral, denoted $\mathcal{H}_\chi$,
\begin{equation}\label{eq:hankelinte}
    \frac{1}{\ep}\int_{\mathcal{H}_{\chi}} \frac{y_0^{(1)}(z)}{\big(\chi_1(z)-w\big)^{\alpha_1}} \e^{-w/\epsilon} \,  \de{w}\sim \sigma_1 y_0^{(1)}(z)\e^{-\chi_1(z)/\epsilon},
\end{equation}
where $\sigma_1=- 2 \pi \mathrm{i}/[\Gamma(\alpha_1)\epsilon^{\alpha_1}]$ is the Stokes multiplier.
Our model problem, with $y_B$ from \eqref{eq:exactBorel}, has $\alpha_1=1/2$ and $y_0^{(1)}(z)=-1/[\sqrt{2}(z-1)]$.
Further details of the above can be found from Crew and Trinh \cite{crew2024resurgent}. 

The contribution from \eqref{eq:hankelinte} is generated across the Stokes line in the physical $\mathbb{C}_z$ plane, denoted $\text{B}>1$, located where $\text{Im}[\chi_1]=0$ and $\text{Re}[\chi_1]>0$. This Stokes line originates at the point where $\chi_1(z_1)=0$; our model problem has $z_1=0$. The smooth interpretation of this transition, across a boundary layer of diminishing width as $\epsilon \to 0$, was derived by formal Borel resummation by Berry \cite{berry_1989}.
It results in an exponentially-small contribution in the transseries of $y(z,\epsilon)$. 

In general, there may be many such contributions analogous to \eqref{eq:hankelinte} to the transseries of $y(z,\epsilon)$, due to multiple singular points of $y_B$. The location of these singularities encode the singulants, $\chi_j(z)$, and the power of the singularity encodes the constants $\alpha_{j}$. Each of these will have corresponding Stokes lines, denoted $\text{B}>j$.
When $y_B$ possesses multiple singular points, the higher-order Stokes phenomenon may emerge, in which these Stokes lines are active in only certain regions of the $\mathbb{C}_z$ plane.
This is discussed next in \S\,\ref{sec:borelHOSP}.

\subsection{Higher-order Stokes phenomenon in the Borel plane}\label{sec:borelHOSP}
We now present a geometrical explanation of HOSP in the Borel plane, with reference to our model example \eqref{eq:themainone}.

The Borel transform, $y_B$ from \eqref{eq:exactBorel}, has two singularities in the Borel plane $\mathbb{C}_w$, a branch point at $w = \chi_1(z)=z^2/2$ and a pole at $w = \chi_2(z)=z-1/2$. In the physical plane $\mathbb{C}_z$, the $\text{B}>1$ Stokes line, located where $\chi_1$ is real and positive, emanates from $z=z_1=0$ at which $\chi_1(z_1)=0$. Similarly, the $\text{B}>2$ Stokes line, located where $\chi_2$ is real and positive, originates from $z=z_2=1/2$ where $\chi_2(z_2)=0$.
Note that the location $z_1=0$ is associated with a singularity in the base expansion of $y(z,\ep)$ from \eqref{eq:base}, but $z_2=1/2$ is not.

\myred{We now examine the Stokes phenomenon produced by the pole at $w=\chi_2$. Note that the branch point at $w=\chi_1$ results in a multi-valued function $y_{B}$ with two sheets. However, the pole at $w=\chi_2$ only exists on one of the two sheets. From the exact Borel solution \eqref{eq:exactBorel}, we have that
\begin{equation}\label{eq:InverseBorelsheet}
y_B(w,z) \sim \left\{\begin{aligned}
&-\frac{1}{(w-\chi_2)}+\cdots\\
&\frac{1}{2(z-1)^2}+\cdots
\end{aligned}\right. \qquad \text{as $w \to \chi_2(z)$}
\end{equation}
when approached on the two corresponding sheets. Thus, when we calculate the inverse Borel transform \eqref{eq:InverseBorel} for the $\text{B}>2$ Stokes line, the result depends on whether the solution with two corresponding behaviours from \eqref{eq:InverseBorelsheet} crosses the integration contour $(0,\infty)$ emanating from the distinguished origin. If the pole crosses the integration contour, the associated residue contribution from the Hankel contour results in a $\text{B}>2$ Stokes phenomenon. Alternatively, if the pole is on the other Riemann sheet, no Stokes phenomenon occurs and we refer to the $\text{B}>2$ Stokes line as being inactive.}

The critical locus where the pole at $w=\chi_2$ is visible in the Borel plane from rays emanating from the distinguished origin, $w = 0$, is the condition that the two singularities are co-linear on a half line:
\begin{equation}\label{eq:higherorderintroline}
  \text{Im}\bigg[ \frac{\chi_2-\chi_1}{\chi_1}\bigg] = 0  \qquad \text{and} \qquad  \text{Re}\bigg[ \frac{\chi_2-\chi_1}{\chi_1}\bigg] \geq 0.
\end{equation}
We define contours that satisfy \eqref{eq:higherorderintroline} as \emph{higher-order Stokes lines}.
Note that these originate in $\mathbb{C}_z$ from the location $z=z_{*}$ at which $\chi_1(z_*) = \chi_2(z_*)$, which is where the $\text{B}>1$ and $\text{B}>2$ Stokes lines intersect. Our model example has $z_*=1$. In the Borel plane this is where the two singularities of $y_B$ are at the same location in $\mathbb{C}_w$ \myred{(though on different Riemann sheets)}.

\begin{figure}
    \centering
    \includegraphics[scale=1]{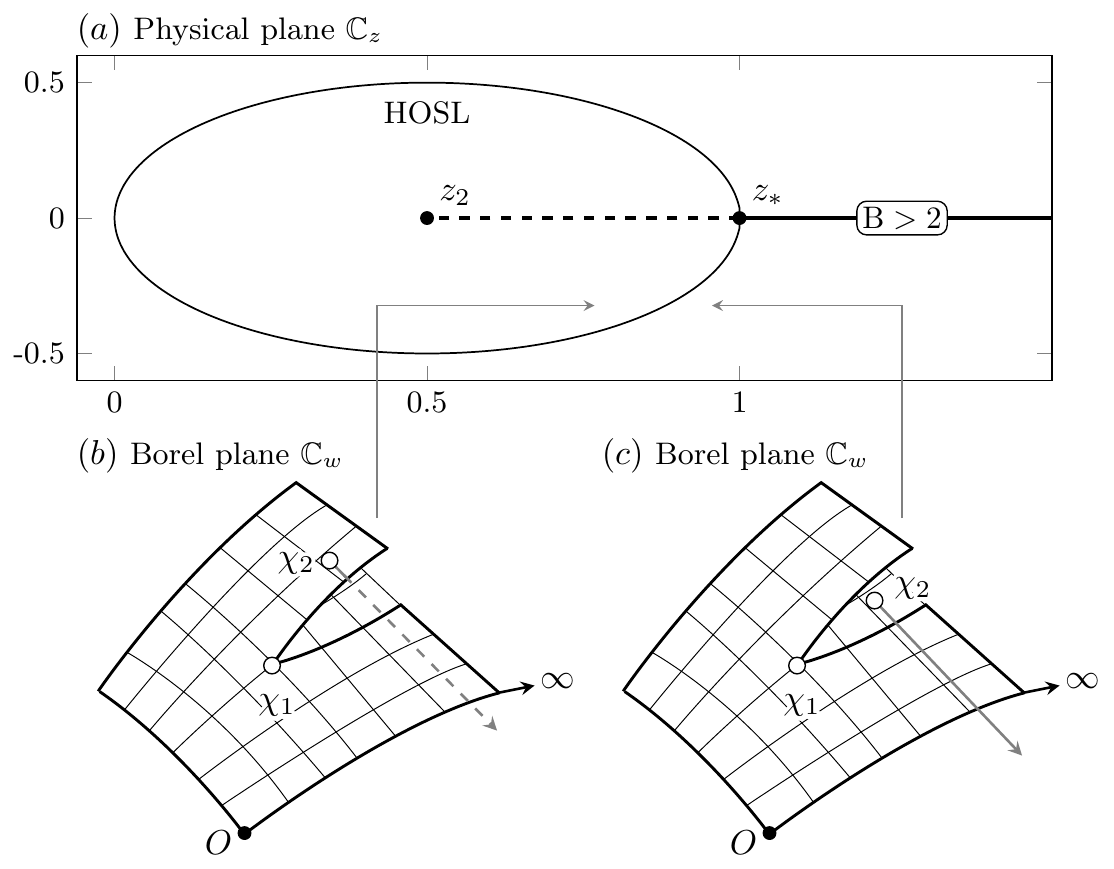}
    \caption{In $(a)$, we plot the geometrical colinearity condition \eqref{eq:higherorderintroline} (HOSL) in $\mathbb{C}_{z}$, which divides the two Borel plane scenarios depicted in $(b)$ and $(c)$.
    In order for the Borel singularity at $w=\chi_2$ to induce the Stokes phenomenon by crossing the integration contour $(0,\infty)$ of the inverse Borel transform \eqref{eq:InverseBorel}, it must be present on the principle Riemann sheet when viewed from the origin, $O$. In ($c$), it remains on the principal sheet and results in a $\text{B}>2$ Stokes line. However, in ($b$), it crosses the real-axis on a higher Riemann sheet and thus no $\text{B}>2$ Stokes phenomenon occurs. }
    \label{fig:BorelIntro}
\end{figure}
The above is a geometrical explanation of HOSP from the perspective of the Borel plane, which produced condition \eqref{eq:higherorderintroline} for a higher-order Stokes line. Thus, it is only in the correct region of the $\mathbb{C}_z$ plane, on one side of \eqref{eq:higherorderintroline}, that the $\text{B}>2$ Stokes line is active. The Stokes line and higher-order Stokes line is illustrated in figure \ref{fig:BorelIntro}. Our explanation here, is similar and complementary to those given by Howls \emph{et al.} \cite{howls_2004} and Chapman and Mortimer \cite{chapman_2005}, in which condition \eqref{eq:higherorderintroline} was also derived.

One of the basic principles of resurgence is that the singularity structure of $y_B$ is encoded in the $n \to \infty$ asymptotics of the late-terms of the expansion $y(z,\epsilon) \sim \sum_{n=0}^{\infty}\epsilon^n y_n(z)$. This follows from Darboux's lemma and extensions thereof. Throughout the rest of this work we shall study the HOSP in $\mathbb{C}_z$ from the perspective of the large-$n$ asymptotics of these late terms, $y_n$. Then, the visibility of $\chi_2$ from the origin in the Borel plane is realised as a Stokes phenomena (in $1/n$) in the expansion of the late terms, which occurs on the colinear locus \eqref{eq:higherorderintroline}. In other words, the geometric condition \eqref{eq:higherorderintroline} will be found to be equivalent to which divergent components are present in the late-term asymptotics of $y_n$.

\section{Asymptotics in the $\mathbb{C}_{z}$ plane}\label{sec:physicalhospsec}
In this section, we perform a direct asymptotic study of the model example \eqref{eq:themainone} in the $\mathbb{C}_{z}$ plane.
We begin in \S\ref{sec:asymptoticinitial} by deriving the early orders of the base expansion, which reveals a singularity at $z=z_1$ only. The divergence of $y_n$ as $n \to \infty$ associated with this singularity, examined in \S\ref{sec:lateanalysisexample}, contains an amplitude function that is singular at the additional location $z=1$. In \S\ref{sec:latelatediv}, we derive lower-order components (in $1/n$) of the late-term divergence, which we denote the late-late-term divergence. This reveals an additional divergent series that induces higher-order Stokes phenomena in the late-terms of the expansion, which produces the anticipated divergent component with $\chi_2=z-1/2$. The smooth behaviour of this is derived in \S\ref{sec:smoothinghosp} through formal arguments.

\subsection{Transseries and divergence}\label{sec:transdiv}
 We begin with a review of the relationship between exponential asymptotic expansions (the {\it transseries}) and the late terms of an asymptotic series. 
 In general, the base expansion \eqref{eq:base} represents only one of many possible contributions to the formal parametric transseries of $y(z; \epsilon)$.
 The full transseries can be written as 
\begin{equation}\label{eq:MainTransseries}
\begin{split}
    y(z;\eps) &= \Big[ y_0(z) + \ep y_1(z) + \cdots + \ep^n y_n(z) + \cdots \Big]   \\
    &\quad+\sum_{i=1
    }\sigma_{i}(z; \ep)\biggl[y^{(i)}_0(z) + \ep y^{(i)}_1(z) + \cdots +\ep^n y^{(i)}_n(z) + \cdots\biggr] \e^{-\chi_i(z)/\ep}.
     \end{split}
\end{equation}
Here, the functions $\chi_i$ are known as the singulants, and $\sigma_i$ are the Stokes multipliers which exhibit the Stokes phenomenon by rapidly changing value across a Stokes line, typically with an error function behaviour \cite{berry_1989}. We will focus on asymptotic series for which $y_n^{(i)}$ diverges factorially as $n \to \infty$. The structure \eqref{eq:MainTransseries} frequently arises in solutions to singularly perturbed ordinary differential equations \cite{dingle_book}.

For the general case above in \eqref{eq:MainTransseries}, the associated Stokes lines of the asymptotic expansion are curves in $\mathbb{C}_{z}$ specified by
\begin{equation}\label{eq:l_chij}
\left\{\begin{aligned}
 \text{B}>j:& \qquad \text{Im}[\chi_j] = 0 \quad \text{and} \quad \text{Re}[\chi_j] \geq 0,\\
   i>j:& \qquad \text{Im}[\chi_j-\chi_i] = 0 \quad \text{and} \quad \text{Re}[\chi_j-\chi_i] \geq 0.
\end{aligned}\right.
\end{equation}
The notation $\text{B}>j$ indicates the Stokes lines induced by the base asymptotic expansion, and $i>j$ indicates those induced by $i$th exponential in \eqref{eq:MainTransseries}. This is the notation used by Howls \emph{et al.} \cite{howls_2004}, and across these curves terms proportional to $\e^{-\chi_j(z)/\ep}$ are switched on in the solution. This behaviour is encoded in the Stokes multipliers, $\sigma_i$.

Let us assume that the components, $y_n(z)$, of the base series are meromorphic functions with singular behaviour at $z=z_{1}$. The leading-order divergence associated with this is given by
\begin{equation}\label{eq:MainDivergence}
    y_{n}(z) \sim  y^{(1)}_0(z)\frac{\Gamma(n+\alpha_1)}{\chi_1(z)^{n+\alpha_1}} \qquad \text{as $n \to \infty$},
\end{equation}
where $\chi_1(z_1)=0$ and typically $\alpha_1$ is constant. We denote \eqref{eq:MainDivergence} as the \emph{late-term divergence}.
Note that the singulant, $\chi_1$, and the amplitude function, $y_0^{(1)}$, in \eqref{eq:MainDivergence} above also appear in one term, $y_0^{(1)}\e^{-\chi_1/\ep}$, of the transseries \eqref{eq:MainTransseries}. This connection between late-term divergence and exponentially subdomiannt terms is fundamental in exponential asymptotics.

Surprisingly, the relationship between the late-terms, $y_n$, of the base expansion, and each of the amplitudes $y_p^{(i)}$ from the $\e^{-\chi_i/\ep}$ transseries component can be found through the consideration of the lower-order components of $y_n$ as $n \to \infty$. Specifically, lower-order corrections to the late-terms of \eqref{eq:MainDivergence} are included with
\begin{equation}\label{eq:factorial123}
    y_{n}(z) \sim \sum_{i=1} \left(\sum_{p=0}^{\infty} y^{(i)}_{p}(z) \frac{\Gamma(n-p+\alpha_i)}{\chi_i(z)^{n-p+\alpha_i}}\right) \qquad \text{as $n \to \infty$}.
\end{equation}
Above, the outer sum indexed by $i$ captures each singulant, $\chi_i$, which contribute to a $\text{B}>j$ Stokes line. The inner sum, indexed by $p$, includes the lower-order correction terms. \myred{This relationship between the $1/n$ late-term expansion and the $\epsilon$ transseries is known as resurgence, to which we refer the reader to the recent review by Aniceto \emph{et al.} \cite{aniceto2019primer}, and references therein. This connection was used by Body \emph{et al.} \cite{body2005exponential} to study HOSP in the physical plane.}
We note that lower-order corrections to the late-term divergence have previously been calculated by {\it e.g.} King and Chapman \cite{king2001asymptotics} but this prior study considered a series in powers of $n^{-1}$, and thus the connection with the transseries \eqref{eq:MainTransseries} was not clear. 

In the following sections, we apply these techniques to our model problem \eqref{eq:themainone}. Here, it is the divergence of $y_p^{(j)}$ in \eqref{eq:factorial123} above as $p \to \infty$ that will lead to the HOSP that was previously derived in \S\ref{sec:borelHOSP} with the Borel plane.

\subsection{Asymptotic expansion}\label{sec:asymptoticinitial}

We begin by seeking the formal transseries expansion for our model ODE \eqref{eq:themainone} by writing $y(z,\epsilon) \sim \sum_{i=0}\e^{-\chi_i(z)/\epsilon}y^{(i)}(z,\epsilon)$,
where $\chi_i(z)$ are the singulant functions, and $y^{(i)}(z,\ep)$ are the corresponding amplitude functions whose expansions were shown in \eqref{eq:MainTransseries}.
Substitution of this ansatz into the differential equation \eqref{eq:themainone} yields three possible solutions for the singulant, which are given by
\begin{equation}\label{eq:chisolutions}
 \chi_0(z)=0, \qquad   \chi_1^{\prime}(z) = z, \qquad \chi^{\prime}_2(z) = 1.
\end{equation}
The first of these, $\chi_0$, corresponds to the base expansion, which will be a standard power series in $\ep$. Furthermore, equations for each of the three amplitude functions, $y^{(i)}(z,\ep)$, are obtained as
\begin{subequations}\label{eq:amplitudeequations}
\begin{align}
\label{eq:amplitudeeq0}
    \ep^2 y^{(0) \prime \prime} + \ep (1+z) y^{(0)\prime} + z y^{(0)} =1,\\
\label{eq:amplitudeeq1}
    \ep y^{(1) \prime \prime} + (1-z) y^{(1)\prime} -y^{(1)}=0,\\
\label{eq:amplitudeeq2}
    \ep y^{(2) \prime \prime} +(z-1) y^{(2)\prime}  =0.
\end{align}
\end{subequations}

We begin by expanding the base solution as $\epsilon \to 0$ by writing $y^{(0)}(z,\ep) \sim \sum_{n=0}^{\infty} \ep^{n} y_{n}(z)$,
for which substitution into equation \eqref{eq:amplitudeeq0} yields the early orders of the solution expansion as
\begin{equation} \label{eq:y0earlyorders}
   y^{(0)}(z,\ep) \sim \frac{1}{z} + \ep\frac{(1+z)}{z^3} + \ep^2 \frac{(2z^2+3z+3)}{z^5}+\cdots.
\end{equation}
Note that there is a singularity at $z=0$, near to which the asymptotic series \eqref{eq:y0earlyorders} reorders. The resultant inner analysis in this region is performed in Appendix \ref{sec:AppInner}, in order to determine unknown constants of our $n \to \infty$ approximation of $y_n$ in \S\ref{sec:lateanalysisexample}. 
At $O(\ep^n)$ in the base equation \eqref{eq:amplitudeeq0}, we find the equation
\begin{equation} \label{eq:y0lateorders}
y_{n-2}^{\prime \prime} +(1+z) y_{n-1}^{\prime} + z y_{n} =0,
\end{equation}
which holds for $n \geq 2$.
As successive orders of the asymptotic expansion are determined by the differentiation of previous orders, the solution $y_n$ will diverge as $n \to \infty$. This divergence is studied in the next section.

\subsection{Late-term divergence}\label{sec:lateanalysisexample}
To capture the factorial-over-power divergence of $y_n$, we consider a late-term ansatz of the form
\begin{equation} \label{eq:fullfactorialpower}
y_n(z) \sim \sum_{p=0}^{\infty} B_p(z)\frac{\Gamma(n-p+\alpha_1)}{\chi_1(z)^{n-p+\alpha_1}} \qquad \text{as $n \to \infty$},
\end{equation}
where $\chi_1$ and $B_p$ are the late-term singulant and amplitude functions, respectively, and $\alpha_1$ is a constant.
In \eqref{eq:fullfactorialpower}, in addition to the leading-order divergence specified by $p=0$, we have also included lower-order components. In fact, there will be additional components of the late-term divergence not captured by the ansatz \eqref{eq:fullfactorialpower}.
These will occur due to the HOSP, in which they are smoothly switched on across a boundary-layer, of diminishing width as $n \to \infty$, that surrounds a HOSL. The key result of this paper is that this phenomenon may be derived through the study of the late-late-terms, $B_p$, as $p \to \infty$, which is performed in \S\ref{sec:smoothinghosp}.

Substitution of the factorial-over-power ansatz \eqref{eq:fullfactorialpower} into the $O(\ep^n)$ equation \eqref{eq:y0lateorders} yields at each order of $\Gamma(n-p+\alpha_1)/\chi_1^{n-p+\alpha_1}$ the equations
\begin{subequations}\label{eq:amplitudeequationsfactorial}
\begin{align}
\label{eq:amplitudeeq0factorial}
p=0: & \qquad (\chi_1^{\prime}-1)(\chi_1^{\prime}-z)=0,\\
\label{eq:amplitudeeq1factorial}
p=1: & \qquad   (1+z-2\chi_1^{\prime})B_0^{\prime}-\chi_1^{\prime \prime}B_0=0,\\
\label{eq:amplitudeeq2factorial}
p \geq 2: & \qquad     B_{p-2}^{\prime \prime} + (1+z-2 \chi_1^{\prime})B_{p-1}^{\prime}-\chi_1^{\prime \prime}B_{p-1}=0.
\end{align}
\end{subequations}

Integration of $\chi_1^{\prime}(z)$ from \eqref{eq:amplitudeeq0factorial} with the matching criteria of $\chi_1(0)=0$ yields two solutions, $\chi_1^{\prime}=z^2/2$ and $\chi_1^{\prime}=z$. However, only the first of these produces an amplitude function that matches to the inner solution from Appendix \ref{sec:AppInnerz0}.
We therefore focus on the first late-term singulant, $\chi_1=z^2/2$, for which the amplitude equation \eqref{eq:amplitudeeq1factorial} may be solved. This yields
\begin{subequations}\label{eq:latetermB0}
\refstepcounter{equation}\label{eq:latetermB0A}
\refstepcounter{equation}\label{eq:latetermB0B}
\begin{equation}
\chi_1(z) =z^2/2 \qquad \text{and} \qquad B_0(z)=\frac{\Lambda_0}{1-z},
  \tag{\ref*{eq:latetermB0A},b}
  \end{equation}
\end{subequations}
where $\Lambda_0$ is a constant of integration. The result of this section is a $\text{B}>1$ Stokes line. Note that solution \eqref{eq:latetermB0B} is singular at $z=1$, and that equation \eqref{eq:amplitudeeq2factorial} for $B_1$ requires differentiation of $B_0$. Hence, $B_1$ will have stronger singular behaviour at $z=1$. This pattern continues and leads to factorial-over-power divergence of $B_p$ as $p \to \infty$, which we denote the late-late-terms of the asymptotic expansion.

\subsection{Late-late-term divergence}\label{sec:latelatediv}
We now determine the divergence of $B_p$ as $p \to \infty$, which arises as a consequence of the singularity at $z=1$ in solution \eqref{eq:latetermB0B} for $B_0$. The power of this singularity will grow as we proceed into the asymptotic series. For instance, equation \eqref{eq:amplitudeeq1factorial} may be solved for $B_1$, yielding
\begin{equation} \label{eq:latetermB1}
B_1(z)=-\frac{\Lambda_0}{(1-z)^3}+\frac{\Lambda_1}{1-z},
\end{equation}
and this growing singular behaviour generates factorial-over-power divergence of $B_p$ as $p \to \infty$.
We begin by considering the factorial-over-power ansatz
\begin{equation} \label{eq:latelatetermBp}
B_p(z) \sim C_0(z)\frac{\Gamma(p+\beta)}{[\tilde{\chi}(z)]^{p+\beta}} \qquad \text{as $p \to \infty$},
\end{equation}
where $\tilde{\chi}$ and $C_0$ are the late-late-term singulant and amplitude functions, and $\beta$ is a constant.
In order to match with the inner solution for $y_n$ near $z=1$ from Appendix \ref{sec:AppInnerz1}, we require $\tilde{\chi}(1)=0$. Substitution of ansatz
\eqref{eq:latelatetermBp} into equation \eqref{eq:amplitudeeq2factorial} yields the equations $\tilde{\chi}^{\prime}=1-z$ and $C_0^{\prime}=0$, which we integrate to find
\begin{equation}\label{eq:latelatesolutions}
\tilde{\chi}(z)=-\frac{(z-1)^2}{2} \qquad \text{and} \qquad C_0(z)=\tilde{\Lambda}.
\end{equation}
Here, $\tilde{\Lambda}=\i/(2\pi)$ is a constant of integration, determined in equation \eqref{eq:latelateconstantapp} of Appendix~\ref{sec:AppInnerz1} through an inner matching procedure at $z=1$.

We now explain how the HOSP is generated from the late-late-terms \eqref{eq:latelatetermBp}.
Due to the divergence of \eqref{eq:latelatetermBp} as $p \to \infty$, optimal truncation of the late-term expansion \eqref{eq:fullfactorialpower} is required at $p=P-1$, where consecutive terms are of the same order. The remainder to this optimally truncated series then displays the Stokes phenomenon across a HOSL, which are contours satisfying the Dingle conditions
\begin{equation}\label{eq:intro-hosl}
 \text{Im}\left[\frac{\tilde{\chi}}{\chi_1}\right] = 0  \qquad \text{and} \qquad \text{Re}\left[\frac{\tilde{\chi}}{\chi_1}\right] \geq 0.
\end{equation}
Thus the geometrical condition \eqref{eq:higherorderintroline} is now realised precisely as a Stokes line, but in the $1/n$ asymptotic series. This transition across the HOSL is smooth, and follows error function dependence over a boundary layer of width $O(n^{-1/2})$. We demonstrate this in \S\ref{sec:smoothinghosp} through Borel resummation of the remainder, which is analogous to the formal method used by Berry \cite{berry_1989} to smooth the regular Stokes phenomenon.

Furthermore, the term switched on across the higher-order Stokes line \eqref{eq:intro-hosl} is another factorial-over-power component of the late-term divergence, with a singulant given by $\chi_1+\tilde{\chi}$. This results in a modified late-term divergence of the base series,
\begin{equation}\label{eq:latetermshospmodified}
    y_{n}(z) \sim \tilde{\sigma}(z,n) C_0(z)\frac{\Gamma(n+\alpha_1+\beta)}{(\chi_1+\tilde{\chi})^{n+\alpha_1+\beta}}+\sum_{p=0}^{P-1} B_{p}(z) \frac{\Gamma(n-p+\alpha_1)}{\chi_1^{n-p+\alpha_1}}.
\end{equation}
In the above, $\tilde{\sigma}$ is the higher-order Stokes multiplier, which rapidly transitions across the HOSL \eqref{eq:intro-hosl} in the smooth manner specified later in equation \eqref{eq:smoothingintegral7}. For our current example, $\tilde{\sigma}=0$ on the left of the HOSL in figure \ref{fig:BorelIntro} and $\tilde{\sigma}=2 \pi \mathrm{i}$ on the right.

Thus, in addition to the $\text{B}>1$ Stokes phenomenon induced by the $\chi_1$ divergence, there will be an additional $\text{B}>2$ Stokes phenomenon induced by
\begin{equation}\label{eq:thenewsingulant}
\chi_2:=\chi_1+\tilde{\chi}=z-1/2,
\end{equation}
which only occurs where $\tilde{\sigma} \neq 0$.
Following Berk \emph{et al.} \cite{berk1982new}, contours that satisfy $\text{Im}[\chi_2]=0$, $\text{Re}[\chi_2] \geq 0$, and $\tilde{\sigma} \neq 0$ are often referred to as {\it new Stokes lines}. The additional $\mathrm{e}^{-\chi_2/\epsilon}$ exponential present in the transseries expansion of $y(z,\epsilon)$, generated by the new Stokes line, is necessary to avoid a contradiction arising from the intersection of regular Stokes lines.

We note that, depending on the specific problem considered, it is also possible for the higher-order Stokes multiplier $\tilde{\sigma}$ to change between two non-zero values due to HOSP, rather than the case above in which one of these values is zero. This corresponds to multiple divergent contributions, each derived from singularities of the early orders, altering one another through HOSP, rather than a new late-term component appearing as in \eqref{eq:latetermshospmodified}.

\subsection{Stokes line structure}
\begin{figure}
    \centering
    \includegraphics[scale=1]{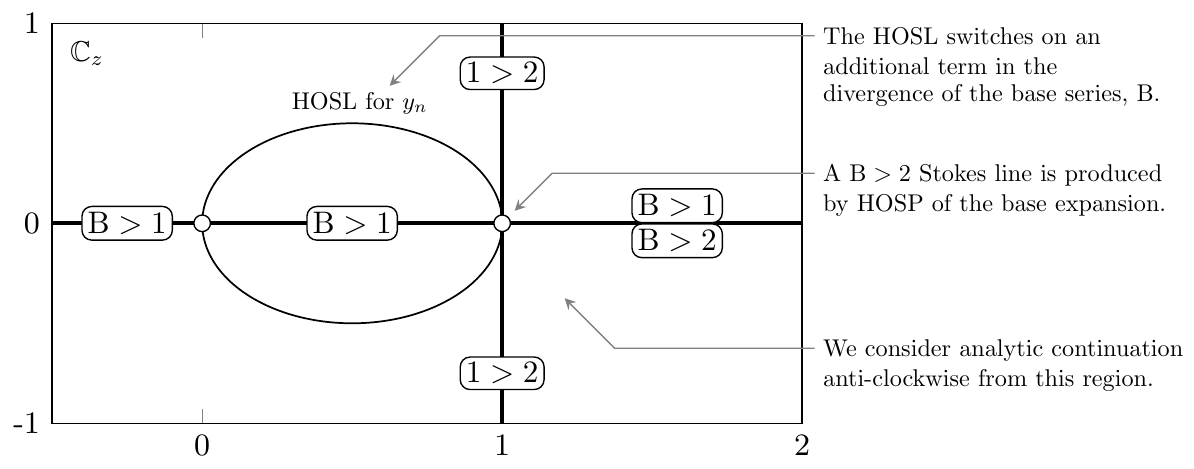}
    \caption{The Stokes line structure for our model problem \eqref{eq:themainone} is shown in $z \in \mathbb{C}$. HOSL are shown with thin lines, and regular Stokes lines with bold lines. \myred{Due to behavioural condition \eqref{eq:themainoneBC}, only the base expansion $\{ \text{B} \}$ is present in the lower-half plane. As we rotate anticlockwise about $z=1$, crossing the $\text{B}>1$ Stokes line ($\text{Re}[z]<1$, $\text{Im}[z]=0$) yields $\{ \text{B},1 \}$ components to the asymptotic transseries in the upper-left quadrant. Crossing $1>2$ ($\text{Re}[z]=1$, $\text{Im}[z]>0$) results in $\{ \text{B} ,1 ,2\}$ present in the upper-right quadrant. Finally, crossing $\text{B}>1$ and $\text{B}>2$ ($\text{Re}[z]>1$, $\text{Im}[z]=0$) returns the asymptotic transseries back to containing only $\{ \text{B} \}$ in the lower-right quadrant.}}
    \label{fig:OurEquation}
\end{figure}

Combined, the Stokes lines for this problem are given by
\begin{equation}\label{eq:allthestokeslines}
\left\{\begin{aligned}
 \text{B}>1:& \qquad  \text{Im}[\chi_1]=0, \quad \text{Re}[\chi_1] \geq 0,\\
  \text{B}>2:& \qquad \text{Im}[\chi_2]=0, \quad \text{Re}[\chi_2] \geq 0,\\
   1>2:& \qquad \text{Im}[\chi_2-\chi_1]=0, \quad  \text{Re}[\chi_2-\chi_1] \geq 0.
\end{aligned}\right.
\end{equation}
The $\text{B}>1$ and $\text{B}>2$ Stokes lines are generated by late-term divergences of the base expansion. The first of these is the factorial-over-power divergence with singulant $\chi_1=z^2/2$ from \eqref{eq:latetermB0A}. The second is the component with $\chi_2=z-1/2$ from \eqref{eq:thenewsingulant} that is switched on across the HOSL. Furthermore, examination of the $O(\e^{-\chi_1/\ep})$ transseries equation \eqref{eq:amplitudeeq1} reveals another divergent series for $y^{(1)}(z,\ep)$. Through the transseries correspondence discussed in \S\ref{sec:transdiv}, this divergent series will have a singulant given by $\tilde{\chi}=\chi_2-\chi_1$, and thus the switching yields a term of $O(\e^{-\chi_2/\ep})$ across a $1>2$ Stokes line, which we included in \eqref{eq:allthestokeslines}. Note that the $\text{B}>2$ Stokes line is active only when the higher-order Stokes multiplier, $\tilde{\sigma}(z,n)$ in \eqref{eq:latetermshospmodified}, is non-zero. These Stokes lines are shown in figure \ref{fig:OurEquation}. For this example, the new $\text{B}>2$ Stokes line generated through HOSP is necessary to avoid a contradiction occurring from rotating about the simple pole at $z=1$.

We have therefore found the following transseries structure with singulants $\chi_1=z^2/2$ and $\chi_2=z-1/2$,
\begin{equation}
\begin{split}\label{eq:transseriesintronewlines}
    y(z,\epsilon) &\sim y_0(z) + \epsilon y_1(z) + \epsilon^2 y_2(z) + \cdots \\
    &\quad+\sigma_1(z,\ep) B_0(z)\e^{-\chi_1/\epsilon} + \sigma_2(z,\ep)C_0(z)\e^{-\chi_2/\epsilon},
    \end{split}
\end{equation}
\myred{where the presence of each component (denoted B, $1$, and $2$) throughout different sectors of the complex-$z$ plane is described in the caption to figure~\ref{fig:OurEquation}.}

\subsection{Higher-order Stokes line smoothing}\label{sec:smoothinghosp}
In this section, we demonstrate that optimal truncation of the late-term ansatz \eqref{eq:fullfactorialpower} yields an exponentially-small (in $n$) remainder that displays the Stokes phenomenon. This remainder corresponds to an additional factorial-over-power component of the late-terms with a singulant given by $\chi_1 + \tilde{\chi}$, as shown previously in \eqref{eq:latetermshospmodified}. Furthermore, we show how this change occurs smoothly across the HOSL \eqref{eq:intro-hosl}, which is located where $\text{Im}[\tilde{\chi}/\chi_1]=0$ and $\text{Re}[\tilde{\chi}/\chi_1] \geq 0$. This analysis closely follows the method used by Berry \cite{berry_1989}, in which the tail of the divergent expansion is re-summed, and examination of the resultant integral reveals error function behaviour localised across the HOSL.

We begin by truncating the late-term ansatz \eqref{eq:fullfactorialpower} at $p=P-1$ by writing
\begin{equation} \label{eq:fullfactorialpowertruncated}
y_n \sim \sum_{p=0}^{P-1} B_p(z)\frac{\Gamma(n-p+\alpha_1)}{\chi_1^{n-p+\alpha_1}} + \underbrace{\sum_{p=P}^{\infty} B_p(z)\frac{\Gamma(n-p+\alpha_1)}{\chi_1^{n-p+\alpha_1}}}_{R_P},
\end{equation}
and consider $P$ to be large. Note that later in equation \eqref{eq:optimalPvalue}, we minimise $R_P$ by specifying $P$ as a function of $n$. In substituting for both the late-late-term divergent ansatz of $B_p$ from \eqref{eq:latelatetermBp} and the integral definition of the Gamma function (since $n-P+\alpha_1$ is also large), we also swap the order of summation and integration. This yields a geometric series, which may be formally resummed to find
\begin{equation}\label{eq:afterresum}
\begin{split}
R_P&=\frac{C_0}{\tilde{\chi}^{\beta}\chi_1^{n+\alpha_1}} \int_0^{\infty} \int_0^{\infty}w^{\beta-1}t^{n+\alpha_1-1}\e^{-w-t} \sum_{p=P}^{\infty}\bigg(\frac{\chi_1 w}{\tilde{\chi} t}\bigg)^p\de{w} \de{t},\\
&=\frac{C_0}{\tilde{\chi}^{\beta+P}\chi_1^{n-P+\alpha_1}}\dashint_0^{\infty} \dashint_0^{\infty}w^{P+\beta-1}t^{n-P+\alpha_1-1}\e^{-w-t} \frac{1}{1-\frac{\chi_1 w}{\tilde{\chi} t}}\de{w} \de{t}.
 \end{split}
\end{equation}
In the above, the principle value integral arises due to the simple pole at $\chi_1 w=\tilde{\chi}t$ introduced by resummation of the infinite geometric series. 
Next, we change integration variables from $(w,t)$ to $(s,f)$ where $s=\chi_1 w /(\tilde{\chi}t)$ and $f=t+w$. This substitution allows for the independent evaluation of the two integrals, one of which is recognised as the definition of the Gamma function. Performing this substitution, and simplifying, yields
\begin{equation}\label{eq:finalintegral}
\begin{split}
R_P&=\frac{C_0}{\chi_1^{n+\alpha_1+\beta}} \dashint_0^{\infty} \frac{s^{P+\beta-1}}{\big(1+\frac{s \tilde{\chi}}{\chi_1}\big)^{n+\alpha_1+\beta}(1-s)} \de{s} \int_0^{\infty}f^{n+\alpha_1+\beta-1} \e^{-f} \de{f}, \\
&=C_0\frac{\Gamma(n+\alpha_1+\beta)}{\chi_1^{n+\alpha_1+\beta}} \dashint_0^{\infty} \frac{s^{P+\beta-1}}{\big(1+\frac{s \tilde{\chi}}{\chi_1}\big)^{n+\alpha_1+\beta}(1-s)} \de{s}.
 \end{split}
\end{equation}
\myred{It is possible to recast the calculation of this section in the hyperterminant language of Olde Daalhuis \cite{daalhuis1996hyperterminants,daalhuis1998hyperterminants} where the integral \eqref{eq:afterresum} may be understood as a level 2 hyperterminant (see equation (3.3) of Olde Daalhuis \cite{daalhuis2009hyperasymptotics}). We refer the reader to recent work by Nemes \cite{nemes2022dingle} and the recent conference talk by Olde Daalhuis \cite{smoothingAdri} for a rigorous treatment of the higher-order Stokes smoothing from this perspective. We now proceed to formally derive the smoothing of \eqref{eq:finalintegral} explicitly for our model example.}

We begin by substituting for $x=s-1$ into \eqref{eq:finalintegral}, which shifts the location of the pole from $s=1$ to $x=0$. As the dominant component of the integrand at $x=0$ lies inside the new range of integration, $-1 \leq x<\infty$, we may also extend the lower range to $-\infty$ without changing the dominant asymptotic behaviour. This gives
\begin{equation}\label{eq:smoothingintegral5}
\begin{split}
R_P&=-C_0\frac{\Gamma(n+\alpha_1+\beta)}{\chi_1^{n+\alpha_1+\beta}} \dashint_{-1}^{\infty} \frac{{(1+x)}^{P+\beta-1}}{x\big[1+(1+x){\tilde{\chi}}/{\chi_1}\big]^{n+\alpha_1+\beta}} \de{x},\\
&\sim -C_0\frac{\Gamma(n+\alpha_1+\beta)}{(\chi_1+\tilde{\chi})^{n+\alpha_1+\beta}} \dashint_{-\infty}^{\infty} \frac{1}{x}\e^{(P+\beta-1)\log{(1+x)}-(n+\alpha_1+\beta)\log{\big[1+\tfrac{x\tilde{\chi}}{\chi_1+\tilde{\chi}}\big]}} \de{x}.
 \end{split}
\end{equation}

We now seek to approximate the integral in \eqref{eq:smoothingintegral5} above, both near the HOSL and close to the pole at $x=0$. In writing
\begin{equation}
 \chi_1 = r \e^{\i \vartheta} \qquad \text{and} \qquad \tilde{\chi} = \tilde{r} \e^{\i \tilde{\vartheta}},
\end{equation}
the HOSL, for which $\tilde{\chi}/\chi_1$ is real and positive, is located at $\tilde{\vartheta}- \vartheta=0$. Thus, in the vicinity of this and the pole at $x=0$, both $x$ and $\tilde{\vartheta}-\vartheta$ are small. We may therefore Taylor expand these components in \eqref{eq:smoothingintegral5} to find 
\begin{equation}\label{eq:smoothingintegralexpansions}
\begin{split}
\e^{(P+\beta-1)\log{(1+x)}} &\sim \e^{(P+\beta-1)\big(x-\tfrac{x^2}{2}+ O(x^3)\big)},\\
\e^{-(n+\alpha_1+\beta)\log{\big[1+\tfrac{x\tilde{\chi}}{\chi_1+\tilde{\chi}}\big]}} &\sim \e^{-(n+\alpha_1+\beta)\big[\tfrac{\tilde{r}x}{r+\tilde{r}} -\tfrac{\tilde{r}^2x^2}{2(r+\tilde{r})^2}+\i\tfrac{r \tilde{r}x(\tilde{\vartheta}-\vartheta)}{(r+\tilde{r})^2}+O(x (\tilde{\vartheta}-\vartheta)^2;x^3)\big]} .
 \end{split}
\end{equation}
The remainder, $R_P$, is minimised when the largest components of \eqref{eq:smoothingintegralexpansions}, $\exp{(Px)}$ and $\exp{(- \frac{\tilde{r}x n}{r+\tilde{r}})}$, cancel upon multiplication. This motivates the specification of our optimal truncation point by
\begin{equation}\label{eq:optimalPvalue}
P=\frac{\tilde{r}n}{r+\tilde{r}}+\tilde{\rho},
\end{equation}
where $0 \leq \tilde{\rho} <1$ ensures that $P$ takes integer values. We may now verify that our prior assumption of $n-P$ being large is self consistent, since $n-P  \sim r n / (r+\tilde{r}) \gg 0$. Substitution of the optimal value of $P$ from \eqref{eq:optimalPvalue} into the integral \eqref{eq:smoothingintegral5}, and using expansions \eqref{eq:smoothingintegralexpansions}, yields real and imaginary components in the integrand. Analogous to the methodology used by Berry \cite{berry_1989}, it may be verified through direct evaluation that the real component of this integral is subdominant to the imaginary component. The leading-order imaginary component then simplifies to give
\begin{equation}\label{eq:smoothingintegral7}
\begin{split}
R_P&\sim \i C_0\frac{\Gamma(n+\alpha_1+\beta)}{(\chi_1+\tilde{\chi})^{n+\alpha_1+\beta}} \dashint_{-\infty}^{\infty} \frac{1}{x}\e^{-\tfrac{r \tilde{r} n x^2}{2(r+\tilde{r})^2}} \sin{\bigg(\frac{r \tilde{r} (\tilde{\vartheta}-\vartheta)n}{(r+\tilde{r})^2}x\bigg)}\de{x},\\
&\sim \sqrt{2 \pi} \i C_0\frac{\Gamma(n+\alpha_1+\beta)}{(\chi_1+\tilde{\chi})^{n+\alpha_1+\beta}} \int_0^{\tfrac{(r \tilde{r})^{1/2}(\tilde{\vartheta}-\vartheta)n^{1/2}}{(r+\tilde{r})}} \e^{-t^2/2}\de{t}.
 \end{split}
\end{equation}

Equation \eqref{eq:smoothingintegral7} is the main result of this section: an explicit solution for the smooth transition of $R_P$ across a boundary layer of width $O(n^{-1/2})$ surrounding the HOSL. Thus, we have that
\begin{equation}\label{eq:smoothingintegral6}
R_P \sim 2 \pi \i C_0(z)\frac{\Gamma(n+\alpha_1+\beta)}{(\chi_1+\tilde{\chi})^{n+\alpha_1+\beta}}
\end{equation}
switches on in the late-terms of the asymptotic solution across the HOSL, which is given by the conditions stated previously in \eqref{eq:intro-hosl}.

It is also interesting to compare our late-late optimal truncation point \eqref{eq:optimalPvalue} to that known in the hyperasymptotic literature for a $1>2$ Stokes phenomenon, \myred{where an exponential of the form $\mathrm{e}^{-\chi_1/\ep}$ is responsible for Stokes switching of an $\mathrm{e}^{-\chi_2/\ep}$ exponential}. From the work of Berry and Howls \cite{berry1991hyperasymptotics}, this value is given by $N_2=r \tilde{r}/(\ep[r+\tilde{r}])$. If we were to calculate the $\text{B}>2$ Stokes phenomenon induced by the new divergent component \eqref{eq:smoothingintegral6} by optimally truncating the base expansion at $N_1 \sim r/\ep$, then \eqref{eq:optimalPvalue} becomes $P\sim r \tilde{r}/(\ep[r+\tilde{r}])$, which matches that known from hyperasymptotics.

\subsection{Summary and relation of asymptotics in $\mathbb{C}_w$ to $\mathbb{C}_z$}
The HOSP may be viewed from the two perspectives we have considered. These are:
\begin{enumerate}[label=(\roman*),leftmargin=*, align = left, labelsep=\parindent, topsep=3pt, itemsep=2pt,itemindent=0pt ]
\item {\it \S\,\ref{sec:modelprob}, The Borel transform.} In the Borel plane, $\mathbb{C}_w$, HOSP occurs when there are at least two singularities, one of which is a branch point with a non-trivial Riemann sheet structure. The HOSL is then a condition for a particular singularity to be `visible' from the origin and therefore cause Stokes phenomena.
\item {\it \S\,\ref{sec:physicalhospsec}, The late-late-term divergence.} Lower-order components of the late-term divergence $y_n$, of $O(\Gamma(n+\alpha_1-p)/\chi_1^{n+\alpha_1-p})$, were considered in \S\ref{sec:latelatediv}. The corresponding amplitude functions diverge as $p \to \infty$; this is referred to as late-late-term divergence. Borel resummation of this optimally truncated expansion, performed in \S\ref{sec:smoothinghosp}, reveals the HOSP.
\end{enumerate}
Singularities of the Borel transform, at $w = \chi_{i}(z)$, correspond to the singulant functions that characterise the divergence of late terms in the transseries.
Moreover, the local behaviour of $y_B$ near such singularities is determined from the constant, $\alpha_{i}$, that appears in divergent representation \eqref{eq:MainDivergence}.
For instance, from the late-term analysis in \S\ref{sec:lateanalysisexample}, the divergence associated with $\chi_1=z^2/2$ has $\alpha_1=1/2$, and that associated with $\chi_2:=\chi_1+\tilde{\chi}=z-1/2$ in \eqref{eq:smoothingintegral6} has $\alpha_2:=\alpha_1+\beta=1$.
Thus, the parametric Borel transform \eqref{eq:exactBorel} had a branch point at $w=\chi_1$ and a pole at $w=\chi_2$.

\section{Further examples of the higher-order Stokes phenomenon}\label{sec:otherexamples}
We now consider three previously studied differential equations whose asymptotic solutions display the HOSP. We demonstrate how this phenomenon may be resolved in these problems through the application of the methods developed in this paper, in which the late-late-term divergence of the asymptotic solution is derived.

The singularly perturbed equations considered in this section are:
\begin{enumerate}[label=(\roman*),leftmargin=*, align = left, labelsep=\parindent, topsep=3pt, itemsep=2pt,itemindent=0pt ]
\item The Pearcey equation in \S\ref{sec:FirstEx}. This is a linear third-order homogeneous equation, previously studied by Howls \emph{et al.} \cite{howls_2004} through consideration of an integral representation of the solution and a hyperterminant analysis.
Our results are summarised in \S\ref{sec:FirstEx} and full details appear in Appendix~\ref{sec:FirstExApp}.
\item A linear second-order inhomogeneous differential equation in \S\ref{sec:SecondEx}. This equation was proposed by Trinh and Chapman \cite{trinh2013new}, who considered an exact series representation for the late-terms of the asymptotic solution.
\item A linear second-order differential equation with an eigenvalue in \S\ref{sec:AppKelvin}, which arises in geophysical fluid dynamics. This equation, and the reusltant eigenvalue divergence, was studied by Shelton \emph{et al.} \cite{shelton2022Hermite,shelton2023kelvin}.
\end{enumerate}
The first two examples contain new Stokes lines which emanate from intersecting Stokes lines. In the third example, the higher-order and regular Stokes lines coincide, which changes the value of the regular Stokes multiplier.
 
\subsection{The Pearcey equation}\label{sec:FirstEx}

The Pearcey function is perhaps the most recognisable subject in connection to the HOSP; it is discussed in the work of Howls \emph{et al.} \cite{howls_2004}, primarily with respect to the integral representation
\begin{equation} \label{eq:Pearceyintegral}
    I(z) = \int_C \e^{-f(s; z)/\ep} \, \de{s} \quad \text{with} \quad f(s; z) = -\im\left(\frac{1}{4}s^4 + \frac{1}{2} s^2 + sz\right).
\end{equation}
Here, $\ep \to 0$ and $C$ is a contour from $\infty \exp(-3\pi\im/8)$ to $\infty \exp(\pi\im/8)$. Alternatively, manipulation of \eqref{eq:Pearceyintegral} yields the linear third-order homogeneous differential equation
\begin{subequations}\label{eq:pearceymainBVP}
\begin{equation}\label{eq:pearceymain}
\ep^3 I^{\prime \prime \prime}(z)-\ep I^{\prime}(z) - \i z I(z)=0,
\end{equation}
\myred{for which we also specify the behavioural condition
\begin{equation}\label{eq:pearceymainBC}
I(z) \sim \frac{1}{z^{1/3}}\exp{\bigg(\frac{3(\sqrt{3}+\mathrm{i})z^{4/3}}{8 \ep}\bigg)} \quad \text{as}~ \lvert z \rvert \to \infty,~ \arg[z]=3 \pi \mathrm{i}/8.
\end{equation}}
\end{subequations}
\myred{A far-field analysis of \eqref{eq:pearceymain} yields the solution $I(z) \sim \alpha z^{-1/3}\mathrm{e}^{3(\mathrm{i}-\sqrt{3})z^{4/3}/(8 \ep)} + \beta z^{-1/3}\mathrm{e}^{-3 \mathrm{i}z^{4/3}/(4 \ep)} + \gamma z^{-1/3}\mathrm{e}^{3(\mathrm{i}+\sqrt{3})z^{4/3}/(8 \ep)}$ as $\lvert z \rvert \to \infty$. In the direction specified by condition \eqref{eq:pearceymainBC}, the third of these components (which appears in the behavioural condition) is subdominant to the other two. Thus, condition \eqref{eq:pearceymainBC} requires $\alpha=0$, $\beta=0$, and $\gamma=1$, and therefore provides three constraints for the third-order ODE.}

The relationship between the integral \eqref{eq:Pearceyintegral} and ODE formulation \eqref{eq:pearceymainBVP} is discussed in section 7.2 of Takei \cite{takei2017wkb}, who considers the BNR equation from Berk \emph{et al.} \cite{berk1982new}, which is obtained by rescaling \eqref{eq:pearceymain}. The general Pearcey integral \eqref{eq:Pearceyintegral} contains a second parameter, assumed to be fixed in the above formulation.

Aided by the integral representation \eqref{eq:Pearceyintegral}, the HOSP can be understood more easily via a steepest descent analysis \cite{howls_2004}. Interestingly, we could not find a direct asymptotic analysis of the HOSP -- as connected to the late terms -- of the Pearcey function in its singularly-perturbed differential equation form. This we do now using the techniques from the previous sections.

Three different WKBJ solutions are permitted to \eqref{eq:pearceymain}, each indexed by a leading exponential, $A_i(z,\ep) \e^{-S_i(z)/\ep}$ for $i = 1, 2, 3$, which yields divergent expansions of the form
\begin{equation}\label{eq:PearceyTransseirs}
\begin{aligned}
I(z) &= \e^{-S_1/\ep}\bigg[A_1^{(0)}+ \cdots + \ep^n \bigg(B_1^{(0)}\frac{\Gamma(n+\alpha)}{(S_2-S_1)^{n+\alpha}} +\tilde{\sigma}_1C_1\frac{\Gamma(n+\alpha)}{(S_3-S_1)^{n+\alpha}}\bigg) + \cdots \bigg]\\
&\phantom{=}+\e^{-S_2/\ep}\bigg[A_2^{(0)}+ \cdots + \ep^n \bigg(B_{2a}^{(0)}\frac{\Gamma(n+\alpha)}{(S_3-S_2)^{n+\alpha}} +B^{(0)}_{2b}\frac{\Gamma(n+\alpha)}{(S_1-S_2)^{n+\alpha}}\bigg) + \cdots \bigg]\\
&\phantom{=}+\e^{-{S_3}/{\ep}}\bigg[A_3^{(0)}+ \cdots + \ep^n \bigg(B_3^{(0)}\frac{\Gamma(n+\alpha)}{(S_2-S_3)^{n+\alpha}} +\tilde{\sigma}_3 C_3\frac{\Gamma(n+\alpha)}{(S_1-S_3)^{n+\alpha}}\bigg) + \cdots \bigg].
\end{aligned}
\end{equation}
\myred{The third term in \eqref{eq:PearceyTransseirs} above corresponds to that specified in behavioural condition \eqref{eq:pearceymainBC}.}
The derivation of each of these components, e.g. {$A_{i}$, $B_{i}$, $C_{i}$}, is given in Appendix~\ref{sec:FirstExApp}. Each of the six factorial-over-power components of \eqref{eq:PearceyTransseirs} induces the Stokes phenomenon on another exponential when the relevant singulant is real and positive. However, two of these divergences are present only when the corresponding higher-order Stokes multiplier, $\tilde{\sigma}_1$ or $\tilde{\sigma}_3$, is non-zero. 
\begin{figure}
    \centering
    \includegraphics[scale=1]{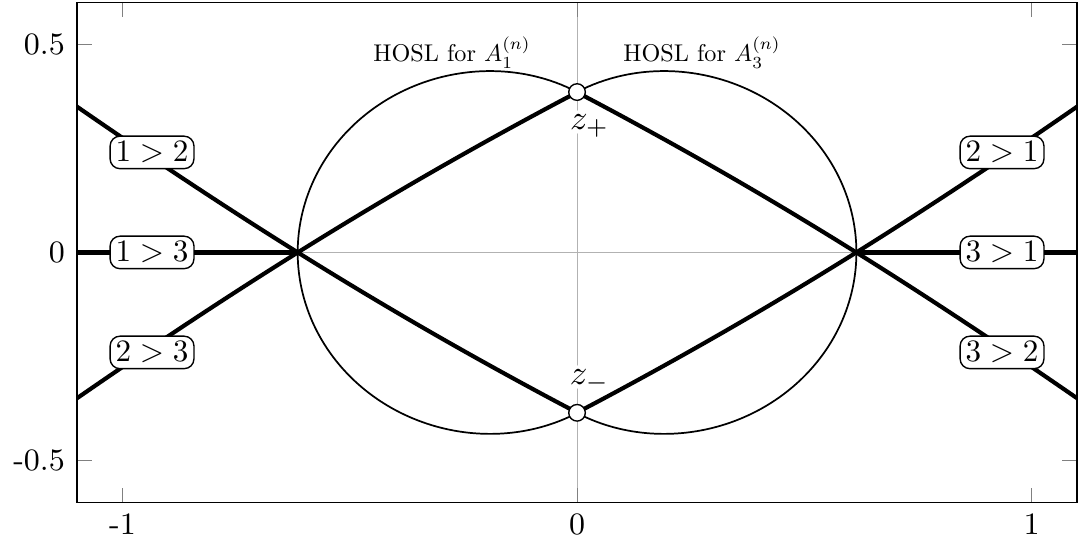}
    \caption{The Stokes line structure is shown in the $z$-plane for the Pearcey equation \eqref{eq:pearceymain}. Higher-order Stokes lines (HOSL), across which the late-terms switch-on a new divergent component, are shown by thin lines. Stokes lines, across which the regular Stokes phenomenon occurs on an exponentially-small component of the asymptotic expansion, are shown as bold lines. On account of the HOSP, two Stokes lines are born (or truncated) at the the point where regular Stokes lines cross. The notation of $i > j$ corresponds to exponential $\e^{S_i/\ep}$ switching on $\e^{S_j/\ep}$ in the asymptotic solution \eqref{eq:PearceyTransseirs}. \myred{We note that this figure is analogous to that presented in figure~4$(a)$ of Howls \emph{et al.} \cite{howls_2004}}.}
    \label{fig:Pearcey}
\end{figure}

The resultant Stokes and higher-order Stokes line structure is shown in figure~\ref{fig:Pearcey}.
Note that there are two singular point at $z_{-}$ and $z_{+}$, which correspond both to singularities in the leading-order asymptotic solution $A_i^{(0)}$, and zeros of the singulant functions, $S_1(z_{-})=S_2(z_-)=S_2(z_+)=S_3(z_{+})=0$. In total, four Stokes lines emanate from $z_{-}$ and $z_{+}$ (shown bold). Stokes lines intersect one another twice, on the real axis. These points also lie on the respective higher-order Stokes lines (shown with thin lines), across which the respective multiplier $\tilde{\sigma}$ switches on to a non-zero value. Thus, beyond the points at which Stokes lines intersect, additional divergences are present in the late-terms of transseries \eqref{eq:PearceyTransseirs}, which leads to two new Stokes lines. These new Stokes lines are marked $1 > 3$ and $3 > 1$ in the figure.
\myred{In the top-right sector, marked in the figure by $(*)$, only the third component, $\{ 3 \}$, of the asymptotic transseries is present, due to behavioural condition \eqref{eq:pearceymainBC}. In rotating clockwise about the point $z=$, we first encounter a $2>1$ Stokes line which has no effect as the second exponential is absent from the transseries. We then cross the $3>1$ Stokes line resulting in the presence of components $\{1,3\}$, the $3>2$ Stokes line giving $\{ 1,2,3\}$, and $2>1$ giving components $\{2,3\}$. Lastly, we return to the original sector by crossing the $3>2$ Stokes line, yielding $\{ 3\}$.}

\subsection{An inhomogeneous second-order equation}\label{sec:SecondEx}
We now study the following second-order singularly perturbed differential equation for $A(z)$,
\begin{subequations}
\begin{align}\label{eq:secondmaineq2}
\ep^2 A^{\prime \prime} + 2 \ep A^{\prime} + (1-z)A=\frac{1}{z-a},\\
\label{eq:secondmaineqBC}
A(z) \to 0 \quad \text{as} \quad \text{Re}[z] \to \pm \infty,~\text{Im}[z]=0.
\end{align}
\end{subequations}
 Equation \eqref{eq:secondmaineq2} was considered as a model problem by Trinh and Chapman \cite{trinh2013new} due to the turning point at $z=1$ and the singularity at $z=a$, and was intended to capture the intersecting Stokes lines appearing in nonlinear gravity-capillary waves travelling on a fluid bounded by lower topography.
In fact, each order of the asymptotic solution to \eqref{eq:secondmaineq2} may be written exactly as a polynomial in powers of $(1-z)^{-1}$ and $(z-a)^{-1}$, for which the constant coefficients are determined as the solution to a recurrence relation. It was conjectured that the analysis of the late-orders behaviour of the recurrence relations could yield details of the HOSP \cite{trinh2013new}, but this was left as an open challenge.

We now apply our general techniques to detect HOSP in this problem. In expanding the solution to equation \eqref{eq:secondmaineq2} as $A(z)=A_0(z) + \ep A_1(z) + \cdots$, the first two orders of this asymptotic expansion are given by
\begin{equation}\label{eq:secondearlysolutions}
A_0=\frac{1}{(1-z)(z-a)} \quad \text{and} \quad A_1= \frac{2}{(1-z)^2(z-a)^2} -\frac{2}{(1-z)^3(z-a)}.
\end{equation}
At $O(\ep^n)$ in equation \eqref{eq:secondmaineq2}, we find $A_{n-2}^{\prime \prime} +2A_{n-1}^{\prime} +(1-z)A_n=0$, which holds for $n \geq 2$.
The late-terms, $A_n$, of the asymptotic expansion will diverge as $n \to \infty$ on account of the singularities at $z=1$ and $z=a$ in \eqref{eq:secondearlysolutions}. 

\subsubsection*{Late-term divergence}
To study the divergence of $A_n$, we consider a factorial-over-power ansatz of the form
\begin{equation}\label{eq:secondexpansion1}
A_n(z) \sim \sum_{p=0}^{\infty} B_p(z) \frac{\Gamma(n-p+\alpha)}{[\chi(z)]^{n-p+\alpha}} \qquad \text{as $n \to \infty$},
\end{equation}
where $\alpha$ is a constant, and $\chi$ and $B_p$ are the singulant and amplitude functions, respectively. Substitution of ansatz \eqref{eq:secondexpansion1} into the $O(\ep^n)$ equation yields at each order of $\Gamma(n-p+\alpha)/\chi^{n-p+\alpha}$ the equations
\begin{subequations}\label{eq:secondeachorderofn}
\begin{align}
\label{eq:secondlateO1}
(\chi^{\prime})^2-2\chi^{\prime}+(1-z)&=0,\\
\label{eq:secondlateOn}
2(1-\chi^{\prime})B_0^{\prime}- \chi^{\prime \prime}B_0&=0,\\
\label{eq:secondlateOnp}
B_{p-2}^{\prime \prime}+2(1-\chi^{\prime})B_{p-1}^{\prime}-\chi^{\prime \prime}B_{p-1}&=0,
\end{align}
\end{subequations}
where \eqref{eq:secondlateOnp} holds for $p \geq 2$.
The quadratic \eqref{eq:secondlateO1} for $\chi^{\prime}$, and \eqref{eq:secondlateOn} for $B_0$ may be solved to find
\begin{equation}\label{eq:secondsingulantsol}
\chi_{\pm}^{\prime}(z) = 1 \pm z^{1/2} \qquad \text{and} \qquad B_0(z) = \frac{\Lambda_0}{(1-\chi_{\pm}^{\prime})^{1/2}},
\end{equation}
where $\pm$ denotes each of the two solutions for $\chi^{\prime}$, and $\Lambda_0$ is a constant of integration.  From the singulant solution in \eqref{eq:secondsingulantsol} above we have $1-\chi^{\prime}_{\pm}=\mp z^{1/2}$, and thus the solution for $B_0$ is singular at $z=0$. This new singularity at $z=0$ will result in the divergence of $B_p$ as $p \to \infty$, which is now studied.

\subsubsection*{Late-late-term divergence}
To capture the divergence of $B_p$ as $p \to \infty$, we consider a typical factorial-over-power ansatz [cf. \eqref{eq:latelatetermBp}]. Substitution of this into equation \eqref{eq:secondlateOnp} yields $\tilde{\chi}^{\prime}_{\pm}(z) = 2(1-\chi_{\pm}^{\prime})$, and hence $\tilde{\chi}_{\pm}(z) = \mp 4 z^{3/2}/3$ upon imposing the condition $\tilde{\chi}_{\pm}(0)=0$. Combined with the solution for the prefactor equation, we find a late-late-term divergence of the form
\begin{equation}\label{eq:secondlatelatediverg}
B_p(z) \sim  \frac{\tilde{\Lambda}}{(\chi^{\prime}_{\pm}-1)^{1/2}} \frac{\Gamma(p+\beta)}{[\mp 4 z^{3/2}/3]^{p+\beta}} \qquad \text{as $p \to \infty$},
\end{equation}
where $\tilde{\Lambda}$ is a constant of integration. This will induce the HOSP in which a new late-term component with singulant $\chi_{\pm} \mp 4 z^{3/2}/3 $ switches on when $\mp 4 z^{3/2}/(3 \chi_{\pm})$ is real and positive.

\subsubsection*{Stokes line structure}
The naive Stokes lines may be determined by integrating the solution $\chi_{\pm}^{\prime}(z) = 1 \pm z^{1/2}$ from \eqref{eq:secondsingulantsol}. As there are two singular points, $z=a$ and $z=1$, that generate divergence, we will have four initial singulants. Only the two generated by the boundary condition $\chi(a)=0$ are now considered to rectify the issue of intersecting Stokes lines found by Trinh and Chapman \cite{trinh2013new}. This yields
\begin{equation}\label{eq:secondchisolutions2}
\chi_1(z) =z+\frac{2}{3}z^{3/2}-a-\frac{2}{3}a^{3/2} \quad \text{and} \quad \chi_2(z)=z-\frac{2}{3}z^{3/2}-a+\frac{2}{3}a^{3/2}.
\end{equation}

Furthermore, the amplitude functions in expansion \eqref{eq:secondexpansion1} for the above two singulants diverge in the manner specified by \eqref{eq:secondlatelatediverg}, on account of a singularity at $z=0$. This new divergent series induces the HOSP in which a new factorial-over-power contribution to $A_n$, with a singulant given by $\chi+\tilde{\chi}$, switches on across the HOSL $\text{Im}[\tilde{\chi}/\chi]=0$ and $\text{Re}[\tilde{\chi}/\chi]>0$. The two new singulants switched on in $A_n$ are therefore
\begin{equation}\label{eq:secondchisolutions3}
\chi_{3}(z)=z-\frac{2}{3}z^{3/2}-a-\frac{2}{3}a^{3/2} \quad \text{and} \quad \chi_{4}(z)= z+\frac{2}{3}z^{3/2}-a+\frac{2}{3}a^{3/2}.
\end{equation}

The new singulants in \eqref{eq:secondchisolutions3}, generated through HOSP, are equal to zero at locations that correspond to regular points of the differential equation \eqref{eq:secondmaineq2}. This was commented upon by Trinh and Chapman \cite[p.~418]{trinh2013new}:
\vspace{1mm}
\begin{quotation}\noindent{\emph{The curiosity, however, is that these singularities do not appear anywhere in the base series, and so they do not seem to be associated with any eventual divergence.}}
\end{quotation}
\vspace{1mm}
As we know now, the divergent contributions of $A_n$ with singulants \eqref{eq:secondchisolutions3} will switch-off across the HOSL; consequently, they do not generate unexpected singular behaviour at the apparent singularity. However, the corresponding singulants still generate divergence in the asymptotic expansions in other regions where they are present.

\begin{figure}
    \centering
    \includegraphics[scale=1]{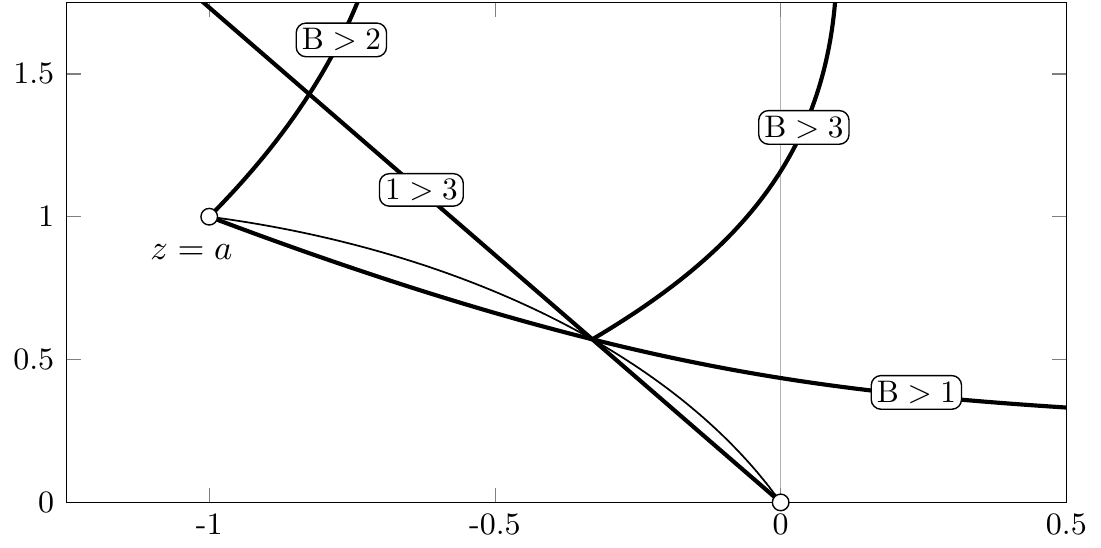}
    \caption{The Stokes line and HOSL structure is shown for $z \in \mathbb{C}$ near the singularity at $z=a$ (circle) for $a=-1+\i$. The turning point at $z=0$ (circle), generates a $1>3$ Stokes line as it is a singularity in the $O(\e^{-\chi_1/\ep})$ expansion. Across the Stokes lines (bold), exponentially-small terms in the asymptotic solution switch on. Across the HOSL (thin line), a new divergent contribution switches on in the late-terms, $A_n$.}
    \label{fig:TrinhChap}
\end{figure}
In figure~\ref{fig:TrinhChap}, we plot the Stokes lines and HOSL associated with the singularity at $z=a$. Two naive Stokes lines, which we denote by $\text{B}>1$ and $\text{B}>2$, are generated by \eqref{eq:secondchisolutions2}. The Airy Stokes line found by evaluating Dingles conditions on $\tilde{\chi}=\mp 4 z^{3/2}/3$, for which $1>3$, is also shown. This results in two points at which Stokes lines intersect. The first, between the $1>3$ and $\text{B}>2$ Stokes lines, is self consistent. However the second, between $1>3$ and $\text{B}>1$, results in a contradiction. The new divergent component with $\chi_3$ generates a $\text{B}>3$ new Stokes line, which resolves this inconsistency.

\subsection{The equatorial Kelvin wave}\label{sec:AppKelvin}

Our last example of the HOSP concerns the study of travelling equatorial Kelvin waves with small latitudinal shear. Here, the system is governed by a second-order singularly perturbed differential equation,
\begin{subequations}
\begin{equation}\label{eq:threeVequation}
\begin{aligned}
\ep^2\bigg[(Y-c-1)(Y-c+1)^2+(Y-c-1)^2(Y-c+1)\bigg] V'' \\
-2Y\bigg[(Y-c-1)(Y-c+1)^2+(Y-c+1)(Y-c-1)^2 +\frac{2 \ep^2(Y-c)^2}{Y} \bigg]V' \\
+\bigg[\big[2Y(Y-c-1)+\ep^2\big](Y-c+1)^2 +(2c-2-\ep^2)(Y-c-1)^2\bigg] V=0,
\end{aligned}
\end{equation}
with boundary conditions
\begin{align}\label{eq:threeVequationBC}
V(0)&= 1,\\
\label{eq:threeVequationDecay}
V(Y)\mathrm{e}^{-Y^2/2\ep^2} &\to 0 \quad \text{as} \quad \text{Re}[Y] \to \pm \infty, ~\text{Im}[Y]=0.
\end{align}
\end{subequations}
The unknown, $V(Y)$, represents the amplitude of a travelling wave ansatz, and $\epsilon$ is the non-dimensional shear. The wavespeed, $c$, is an eigenvalue whose value is obtained from enforcing boundary condition \eqref{eq:threeVequationBC}. \myred{The decay condition \eqref{eq:threeVequationDecay} demands that only the base expansion is present in the far-field of the real axis.} For further details of this model, see Griffiths \cite{griffiths2008limiting} and Shelton \emph{et al.} \cite{shelton2023kelvin}. We consider a solution expansion of the form $V(Y)=V_0 + \ep^2 V_1 + \ep^4 V_2+\cdots$. Enforcing the boundary condition at each order of $\ep$ is only possible if the eigenvalue is also expanded as $c=c_0 + \ep^2 c_1 + \ep^4 c_2+\cdots$.

At leading order, once the boundary condition $V_0(0)=1$ is satisfied, the solutions are given by
\begin{equation}\label{eq:threeO1soleig}
V_0(Y)= (1-Y)^{\frac{1}{2}}\e^{Y/2}  \qquad \text{and} \qquad c_0=1.
\end{equation}
This singularity at $Y=1$ results in a divergent series expansion. Furthermore, the eigenvalue expansion, $c=1+\frac{1}{2}\ep^2 - \frac{1}{8}\ep^4 + \cdots$, will also diverge.

\subsubsection{Late-term divergence}
Due to the nonlinearity between the solution, $V$, and eigenvalue, $c$, in equation \eqref{eq:threeVequation}, the $O(\ep^{2n})$ equation will contain a countably infinite number of terms as $n \to \infty$, such as homogeneous terms, for instance $V_n$ and $V_{n-1}^{\prime}$, and inhomogenous terms such as $c_n$. Both of these contribute to the leading-order divergent solution, 
\begin{equation}\label{eq:A0sol}
V_n(Y) \sim B_0^{(c_n)}(Y)\Gamma(n+\alpha)+ B_0(Y)\frac{\Gamma(n+\gamma)}{(1-Y^2)^{n+\gamma}},
\end{equation}
which contains one term associated with the divergence of $c_n$, and another generated by the singular behaviour at $Y=1$. In \eqref{eq:A0sol}, $\alpha$ and $\gamma$ are constants, and the singulant function $\chi=1-Y^2$ satisfies the condition $\chi(1)=0$. 

Similar to the example from \S\ref{sec:SecondEx}, there is an apparent paradox in that a singularity is seemingly present at $Y=-1$ in the factorial-over-power representation \eqref{eq:A0sol}, which is absent from the early orders of the expansion. This may be resolved by considering lower-order components of $B_0^{(c_n)}$, which induce the HOSP in which the second component of \eqref{eq:A0sol} is switched off before $Y=-1$ is reached. 

\subsubsection{Late-late-term divergence}
In the study by Shelton \emph{et al.} \cite{shelton2023kelvin} the consequences of the HOSP in this problem, an inactive Stokes line and also an atypical Stokes multiplier, were assumed. We now demonstrate how this may be derived from the application of the methods developed in this present paper. 

\begin{figure}
    \centering
    \includegraphics[scale=1]{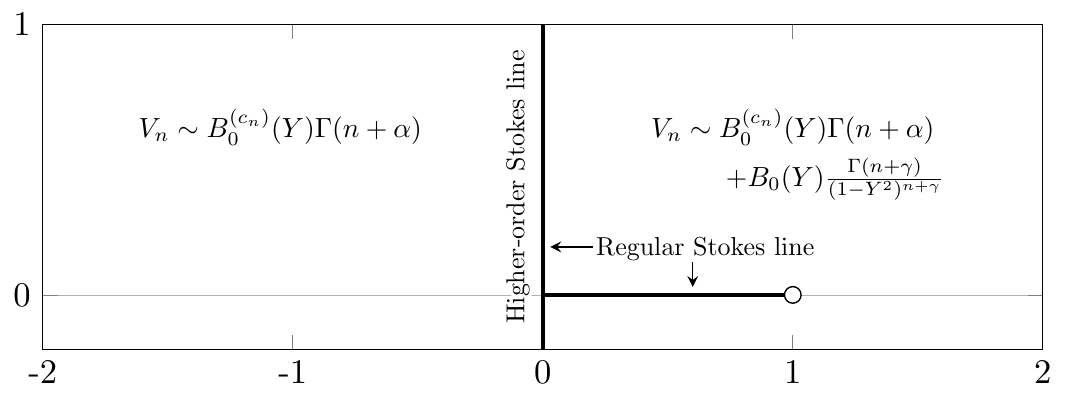}
    \caption{The Stokes line structure is shown for the Kelvin wave equation \eqref{eq:threeVequation}. The regular Stokes line (bold) for $\chi=1-Y^2$ lies along the imaginary axis and part of the real-axis. There is also a HOSL along the imaginary axis, which causes the regular Stokes line on the negative real axis to be switched off.}
    \label{fig:Eigenvalue}
\end{figure}

The late-terms of the eigenvalue expansion will diverge in the factorial-over-power manner of $c_n \sim \delta\Gamma(n+\alpha)$, where $\delta$ and $\alpha$ are constants. The particular solution generated by this eigenvalue divergence can then be determined with a late-term ansatz of the form
\begin{equation}\label{eq:thirdexpansion2}
V_n(Y) \sim \sum_{p=0}^{\infty} B^{(c_n)}_p(Y) {\Gamma(n-p+\alpha)} \qquad \text{as $n \to \infty$}.
\end{equation}
It may be determined that $B_0^{(c_n)} \sim 1/Y$ as $Y \to 0$, which leads to late-late-term divergence of the form
\begin{equation}\label{eq:thirdlatelatediverg}
B_p^{(c_n)}(Y) \sim C(Y) \frac{\Gamma(p+\beta)}{[\tilde{\chi}(Y)]^{p+\beta}} \qquad \text{as $p \to \infty$}.
\end{equation}
The requirement that $\tilde{\chi}(0)=0$ then yields $\tilde{\chi}(Y)=-Y^2$. Thus, across the HOSL $\text{Im}[-Y^2/1]=0$ and $\text{Re}[-Y^2/1]>0$, a new factorial-over-power component of $V_n$ with a singulant given by $\chi=1-Y^2$ will switch on. This is the same singulant as that in \eqref{eq:A0sol} generated by the singularity at $Y=1$. Thus, the apparent contradiction of $\chi=1-Y^2$ containing an undesired singularity at $Y=-1$ is resolved, since the prior expression in \eqref{eq:A0sol} switches off as we cross the HOSL. This is shown in figure~\ref{fig:Eigenvalue}.

\section{Discussion}\label{sec:Discussion}

Although somewhat obscure, the higher-order Stokes phenomenon is generic in the asymptotics of singular perturbation problems. As we have noted, it accompanies situations where the late-term approximation of an $\ep$-expansion, say $y_n$, contains an additional divergent series in inverse powers of $n$. Consequently Stokes phenomena occurs (in $n$), switching on or off late terms and hence Stokes lines themselves. Alternatively, the genericity of HOSP can be interpreted as a consequence of the Borel plane possessing a branched structure. The HOSP often arises in situations where:  (i) Stokes lines cross; (ii) additional singularities, not present in early orders, are predicted in the late terms; or (iii) a naive divergent series is unable to satisfy necessary boundary conditions.

The methods developed in this paper have been presented for linear ODEs. However, they are also applicable to more general nonlinear problems and PDEs. Our approach is similar, and in some cases more general, than those of Chapman and Mortimer \cite{chapman_2005}, Body \emph{et al.} \cite{body2005exponential}, and Howls \emph{et al.} \cite{howls_2004}; in ours, we have focused on the analysis of the $1/n$ divergent series of the late terms and its subsequent Stokes phenomena. Additionally, we have shown that the HOSP and second-generation Stokes phenomenon are both caused by the same type of divergence, either in the late-late-terms, or found in the subsequent transseries expansion. Thus, any HOSP that \emph{e.g.} generates a $\text{B}>2$ Stokes line must necessarily accompany a second-generation $1>2$ Stokes line.

One of the main complications of nonlinearity is the matching procedure to an inner solution that determines components of the late- and late-late-term divergences. For our linear examples, the series expansion for the outer limit of the inner solution could be determined explicitly. This will not be the case in general for nonlinear problems, where the associated recurrence relation is often evaluated numerically. In the study of PDEs, integration of the singulant and amplitude equations from the late- and late-late-term ansatzes can also be challenging \cite{chapman_2005, body2005exponential}. The singularity structure of in the Borel plane will also be more complicated, for which additional Riemann sheets are expected.

\myred{The singularity structure in the Borel plane can be explored for nonlinear problems via numerical analytic continuation techniques such as Pad\'e approximants \cite{costin2022uniformization}, a method used to analytically continue truncated power series. This can be used as a useful heuristic to understand the parametric Borel plane singularity structure for a given parametric series. It involves computing $2N+1$ perturbative terms of the truncated power series $y^{2N}_B(w;z) = y_0(z) + y_1(z)w + \cdots + y_{2N}(z)w^{2N}/(2N)!$.
We then compute the off-diagonal $[N-1:N]$ Pad\'e approximant for the fixed values of $z \in \mathbb{C}_z$ of interest. That is, we write $P_{[N-1:N]}(w;z) = P_{N-1}(w;z)/Q_{N}(w;z)$, where $P_{N-1}$ and $Q_{N}$ are polynomials in $w$, whose coefficients depend on $z$, chosen to agree with the truncated power series at a given order.
\begin{figure}
    \centering
    \includegraphics[scale=1]{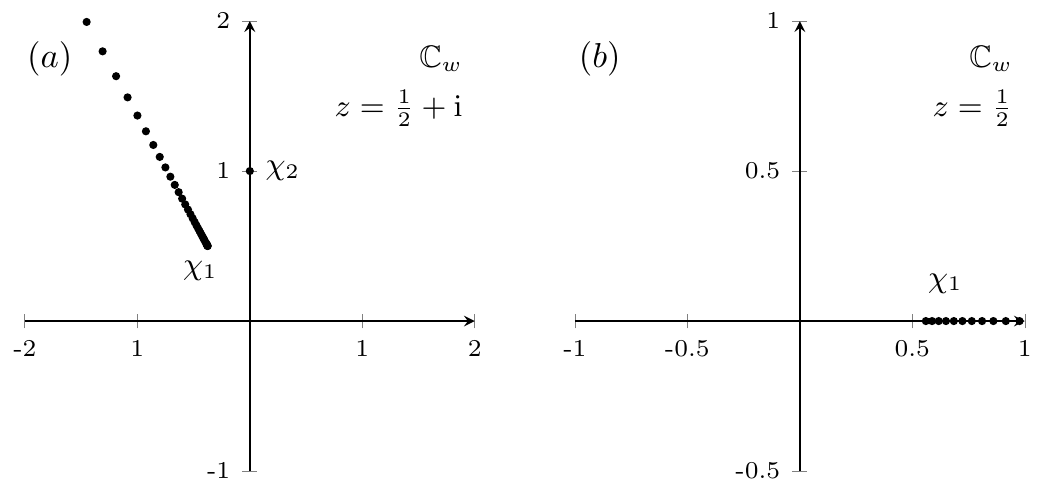}
    \caption{The Borel plane, $\mathbb{C}_w$, is shown for our minimal example \eqref{eq:themainone_intro} in ($a$) for $z=1/2+\i$, and in ($b$) for $z=1/2$. The singularities (nodes) correspond to $w = \chi_1(z) = z^2/2$ and $w = \chi_3(z) = z - 1/2$, and are computed from a Pad\'{e} approximant with $N = 250$. \myred{Note the characteristic coalescing poles that the rational approximation assigns to branch points. The situation $(a)$ corresponds to the region `above' the HOSL in figure~\ref{fig:OurEquation} while $(b)$ corresponds to `below' the HOSL.}}
    \label{fig:BorelExp}
\end{figure}
  We demonstrate this in figure \ref{fig:BorelExp} for our model equation \eqref{eq:themainone_intro}. The interesting region to explore is near the HOSL where one may see poles disappearing behind branch cuts, realising the HOSP. We plot the Borel plane inside the HOSL for $z=1/2 + \i$ in \ref{fig:BorelExp}$(a)$, and outside for $z=1/2$ in \ref{fig:BorelExp}$(b)$. One interesting fact to note is that in $(b)$ when $z=1/2$, the Borel singularity $w=\chi_2(1/2)=0$; however the series $y_B$ is convergent here as the singularity $w=0$ lies on a higher sheet above the distinguished origin (cf. the exact solution \eqref{eq:exactBorel}).
  }

In addition to the development of methodology, another contribution of our work has been to highlight the minimal example \eqref{eq:themainone_intro}. Similar to the example given in Olde Daalhuis \cite{daalhuis2004higher}, such cases are reminders of the subtle and complicated nature of divergence and transseries, even in apparently benign problems. There are a number of exciting ongoing works related to HOSP that have emerged from a recent Isaac Newton Institute programme, notably by Olde Daalhuis \cite{smoothingAdri} on the smoothing of Stokes multipliers of an $\ep$ expansion, King \cite{iniKing} on coincident Stokes lines and HOSP, and Nemes \cite{nemes2022dingle} on atypical Stokes multipliers due to HOSP. This is an exciting and advanced area of research with still more to discover.

\mbox{}\par
{\bf \noindent Acknowledgements}.
We thank the Isaac Newton Institute for Mathematical Sciences, Cambridge, for the programme ``Applicable resurgent asymptotics: towards a universal theory", during which many insightful discussions took place. In particular, we are thankful for discussions with S. J. Chapman (Oxford), G. Nemes (Tokyo Metropolitan University), A. Olde Daalhuis (Edinburgh), and many other participants of the programme. This research was conducted while visiting the Okinawa Institute of Science and Technology (OIST) through the Theoretical Sciences Visiting Program (TSVP). We would also like to thank the anonymous referees for their detailed comments on the contents and structure of our manuscript.

\providecommand{\noopsort}[1]{}

\appendix

\section{Inner analyses near $z=0$ and $z=1$}\label{sec:AppInner}
The purpose of this section is to determine the unknown constants appearing in the late-terms, $y_n$, of the asymptotic solution of \S\ref{sec:modelprob}. We note that while each function in \eqref{eq:fullfactorialpower}, $B_p$ for $p \geq 0$, will have an unknown constant of integration, $\Lambda_p$, not all of these are required. Only knowledge of $\Lambda_0$ and the constant appearing in the divergent form of $B_p$ as $p \to \infty$, $\tilde{\Lambda}$ from \eqref{eq:latelatesolutions}, are required.

\subsection{Inner analysis near $z=0$}\label{sec:AppInnerz0}
We begin by deriving the inner equation at $z=0$, near to which the early orders of expansion \eqref{eq:y0earlyorders} reorder. This is due to a turning point in the original differential equation that generates a singularity in the leading-order solution. Since $y_0 \sim z^{-1}$ and $\ep y_1 \sim \ep z^{-3}$ as $z \to 0$ from equation \eqref{eq:y0earlyorders}, the early orders of the base expansion reorder when $z \sim \ep^{1/2}$. To study the solution in this region, we introduce the inner variable, $\hat{z}$, and the inner solution, $\hat{y}$, by
\begin{subequations} \label{eq:AppInnerVar}
\refstepcounter{equation}\latexlabel{eq:AppInnerVarA}
\refstepcounter{equation}\latexlabel{eq:AppInnerVarB}
\begin{equation}
  z= \ep^{1/2} \hat{z} \qquad \text{and} \qquad   y_{\text{outer}}= \frac{1}{ (\ep^{1/2} \hat{z})} \hat{y}_{\text{inner}}(\hat{z}).
\tag{\ref*{eq:AppInnerVarA},b}
\end{equation}
\end{subequations}
The inner region is defined by $\hat{z}=O(1)$, and our definition of the inner solution in \eqref{eq:AppInnerVarB} ensures that $\hat{y}_{\text{inner}} =O(1)$ in the inner region. \myred{Note that while $z=1/2+O(\ep)$ produces a valid dominant balance in equation \eqref{eq:amplitudeeq0}, corresponding to the base expansion reordering with $\chi_2$ exponential $\mathrm{e}^{-(z-1/2)/\ep}$, this latter term is not present near $z=1/2$ on account of the higher-order Stokes phenomenon.} The leading-order inner equation may then be derived by substituting equations \eqref{eq:AppInnerVar} into the outer equation \eqref{eq:amplitudeeq0}, yielding
\begin{equation}\label{eq:AppInnerEq}
 \frac{1}{\hat{z}} \dd{\hat{y}}{\hat{z}} +\bigg(1-\frac{1}{\hat{z}^2} \bigg) \hat{y}= 1.
\end{equation}

\myred{The inner-limit of the early-orders of the outer expansion from \eqref{eq:y0earlyorders} is given by $\hat{y} \sim  1+{\hat{z}^{-2}} + 3{\hat{z}^{-4}}+\cdots$. Knowledge of the inner solution of equation \eqref{eq:AppInnerEq} is desired, in order to match with the late-terms of the outer expansion and determine the unknown constant of integration. Equation \eqref{eq:AppInnerEq} has the exact solution $\hat{y}(\hat{z})=\hat{z}\mathrm{e}^{-\hat{z}^2/2}[c_1-\mathrm{i}\sqrt{\pi/2}\text{erf}(\mathrm{i}\hat{z}/\sqrt{2})]$, for which expansion as $\hat{z} \to \infty$ yields at algebraic orders
\begin{equation}\label{eq:AppInnerSol3}
\hat{y}(\hat{z}) \sim 1+\sum_{n=1}^{\infty} \frac{\Gamma(2n)}{2^{n-1}\Gamma(n)}\frac{1}{\hat{z}^{2n}}.
\end{equation}

Typically for nonlinear problems, in which exact solutions are absent, a series expansion of the form \eqref{eq:AppInnerSol3} must be assumed. The coefficients of this would be determined numerically via a recurrence relation to a reasonably large index ($n \approx 100$), for which matching with the outer solution yields an approximate value for the constants of integration.
}

\myred{We now match with the outer divergent solution by writing the inner solution \eqref{eq:AppInnerSol3} in terms of the outer variable, $z$, and outer solution, $y$, from equation \eqref{eq:AppInnerVar}. This outer-limit of the inner solution \eqref{eq:AppInnerSol3} is given by $y\sim  \sum_{n=0}^{\infty} \ep^n \Gamma(2n)2^{1-n} z^{-(1+2n)} /\Gamma(n)$.
The inner limit (as $z \to 0$) of the outer divergent solution \eqref{eq:fullfactorialpower}, with amplitude function \eqref{eq:latetermB0} is given by $y_n \sim  \Lambda_0 \Gamma(n+\alpha)(z^2/2)^{-n-\alpha}$.}
Matching of the $O(\ep^n)$ component requires $\alpha=1/2$, which then gives an expression for the constant, $\Lambda_0$. Thus, we have found
\begin{equation}\label{eq:appendixalpha}
\alpha=1/2 \qquad \text{and} \qquad \Lambda_0 = \lim_{n \to \infty}\bigg( \frac{\Gamma(2n)}{2^{2n-1/2}\Gamma(n)\Gamma(n+1/2)}\bigg)= \frac{1}{\sqrt{2 \pi}}.
\end{equation}

\subsection{Inner analysis near $z=1$}\label{sec:AppInnerz1}
The constants $\beta$ and $\tilde{\Lambda}$, from the divergence of $B_p$ as $p \to \infty$ in \eqref{eq:latelatetermBp} are now determined by matching with an inner solution at $z=1$. One method is to consider a boundary layer of diminishing width as $n \to \infty$ on account of the late-term reordering $B_0 \Gamma(n+\alpha)/\chi^{n+\alpha} \sim B_1 \Gamma(n+\alpha-1) / \chi^{n+\alpha-1}$. However, we note that a power series expansion,
\begin{equation}\label{eq:newtransseriesexp}
y^{(1)}(z,\ep) = \sum_{p=0}^{\infty} \ep^p B_p(z),
\end{equation}
in the $\chi=z^2/2$ transseries equation \eqref{eq:amplitudeeq1} yields at each order of $\ep$ the same equations as those satisfied by $B_p$ in \eqref{eq:amplitudeeq1factorial} and \eqref{eq:amplitudeeq2factorial}. Thus, we may determine the unknown constants in the divergent form of $B_p$ through the consideration of a boundary layer, of diminishing width as $\ep \to 0$, at $z=1$ in the governing equation \eqref{eq:amplitudeeq1} for $y^{(1)}$.

Since the early orders of expansion \eqref{eq:newtransseriesexp} reorder when $\Lambda_0 (1-z)^{-1} \sim -\ep \Lambda_0 (1-z)^{-3}$, and $y^{(1)} \sim \Lambda_0 / (1-z)$, we introduce the inner variable, $\hat{z}$ and inner solution, $\hat{y}$, via
\begin{subequations}\label{eq:newtranseriesinnerthings}
\refstepcounter{equation}\latexlabel{eq:newtranseriesinnerthingsA}
\refstepcounter{equation}\latexlabel{eq:newtranseriesinnerthingsB}
\begin{equation}
(1-z) = \ep^{1/2} \hat{z} \qquad \text{and} \qquad y^{(1)}=\frac{\Lambda_0}{\ep^{1/2}\hat{z}}\hat{y}(\hat{z}).
\tag{\ref*{eq:newtranseriesinnerthingsA},b}
\end{equation}
\end{subequations}

The leading-order inner equation for $\hat{y}$ may be found by substituting relations \eqref{eq:newtranseriesinnerthings} into the transseries equation \eqref{eq:amplitudeeq1} for $y^{(1)}$, and retaining the dominant terms as $\ep \to 0$. This yields
\begin{equation}\label{eq:AppInnerEqz1}
\frac{1}{\hat{z}}\dd{^2\hat{y}}{\hat{z}^2}-\bigg(1+\frac{2}{\hat{z}^2}\bigg)\dd{\hat{y}}{\hat{z}}+\frac{2}{\hat{z}^3}\hat{y}=0.
\end{equation}
\myred{This has the exact solution $\hat{y}(\hat{z})=\hat{z}\mathrm{e}^{\hat{z}^2/2}[c_2-c_1\sqrt{\pi/2}\text{erf}(\hat{z}/\sqrt{2})]$, for which expansion as $\hat{z} \to \infty$ and using the matching condition $c_1=1$ yields at algebraic orders
\begin{equation}\label{eq:AppInnerEqz1_again}
\hat{y}(\hat{z}) \sim \sum_{p=0}^{\infty} \frac{(-1)^p\Gamma(2p)}{2^{p-1}\Gamma(p)} \hat{z}^{-{2p}}.
\end{equation}

Now that the outer-limit of the inner solution near $z=1$ is known from \eqref{eq:AppInnerEqz1_again}, we may match it to the inner-limit of the outer solution, in terms of the outer variable, $z$. This inner limit of the late-late terms \eqref{eq:latelatetermBp} is given by $B_p \sim \tilde{\Lambda}(-1)^{p+\beta}2^{p+\beta}\Gamma(p+\beta)(1-z)^{-2p-2\beta}$.} Comparison to the $O(\ep^{p})$ component of \eqref{eq:AppInnerEqz1_again}, when written in outer variables, then yields the constants $\beta$ and $\tilde{\Lambda}$ as
\begin{equation}\label{eq:latelateconstantapp}
\beta=1/2 \qquad \text{and} \qquad \tilde{\Lambda}= \frac{\i}{2 \pi}.
\end{equation}

\section{Details for the Pearcey equation}\label{sec:FirstExApp}
Equation \eqref{eq:pearceymain} permits three WKBJ solutions of the form $A_m(z,\ep) \exp{(-S_{m}(z) / \ep)}$, where $m=1$, $2$, and $3$. Substitution of this ansatz into \eqref{eq:pearceymain} yields the following equations for the singulants, $S_m^{\prime}$, and amplitude function, $A_m$:
\begin{subequations}\label{eq:pearcyafterWKB}
\begin{align}
\label{eq:pearceyexpsol}
\big[S_m'(z)\big]^3-S_m'(z)+\i z=0, \\
\label{eq:pearceyampeq}
\ep^2 A_m^{\prime \prime \prime}-3 \ep S_m^{\prime} A_m^{\prime \prime}+\Big(3(S_m^{\prime})^2-3 \ep S_m^{\prime \prime}-1\Big)A_m^{\prime}+\Big(3S_m^{\prime}S_m^{\prime \prime}-\epsilon S_m^{\prime \prime \prime}\Big)A_m=0.
\end{align}
\end{subequations}

It will be seen that the asymptotic expansions for the amplitude functions, $A_m$, diverge. The Stokes phenomenon that occurs as a consequence of this divergence will switch on one of the other initial singulants. For instance, across a $1>2$ Stokes line, we will have $ A_1 \exp{(-S_1/\ep)} \mapsto A_1 \exp{(-S_1/\ep)} + A_2 \exp{(-S_2/\ep)}$. However, this $1>2$ Stokes line will intersect a $2>3$ Stokes line, for which a new Stokes line is required in order to avoid a contradiction. We now demonstrate how this new Stokes line, and the HOSL that generates it, may be derived through determination of the late-late-term divergence.

In expanding $A_m = A_m^{(0)} + \ep A_m^{(1)} + \cdots$ in equation \eqref{eq:pearceyampeq}, we find the following equations at each order of $\ep$:
\begin{subequations}\label{eq:pearcyeachorderofamp}
\begin{align}
\label{eq:pearceyexpeqO1}
O(1):& \qquad \mathcal{D}\big[A_m^{(0)}\big] \equiv \Big[3(S_m^{\prime})^2-1\Big]A^{(0)\prime}_m+3 S_m^{\prime} S_m^{\prime \prime} A^{(0)}_m=0,\\
\label{eq:pearceyexpeqOep}
O(\ep):& \qquad \mathcal{D}\big[A_{m}^{(1)}\big] =3 S_m^{\prime}A_m^{(0)\prime \prime}+3 S_m^{\prime \prime}A_m^{(0)\prime}+S_m^{\prime \prime \prime}A_m^{(0)},\\
\label{eq:pearceyexpeqOepn}
O(\ep^n):& \qquad \mathcal{D}\big[A_{m}^{(n)}\big] =3 S_m^{\prime}A_{m}^{(n-1)\prime \prime}+3 S_m^{\prime \prime}A_{m}^{(n-1)\prime}+S_m^{\prime \prime \prime}A_{m}^{(n-1)}-A_{m}^{(n-2)\prime \prime \prime},
\end{align}
\end{subequations}
where $\mathcal{D}$ is the linear operator defined in the leading-order equation \eqref{eq:pearceyexpeqO1}.
The leading-order solutions is given by
\begin{equation}\label{eq:pearceyO1sola}
A^{(0)}_m(z)=\frac{C^{(0)}_m}{(3[S_m^{\prime}(z)]^2-1)^{1/2}},
\end{equation}
where $C_m^{(0)}$ is a constant of integration. The singular behaviour of \eqref{eq:pearceyO1sola} may be determined by solving the cubic equation \eqref{eq:pearceyexpsol} for $S_m^{\prime}$, which yields the three solutions
\begin{equation}\label{eq:valuesofS123}
S^{\prime}_1(z)= \frac{12+p^2(z)}{6p(z)} \quad \text{and} \quad S^{\prime}_{2,3}(z)=-\frac{(1 +a \i \sqrt{3})p^2(z)+12(1 -a \i \sqrt{3})}{12p(z)}.
\end{equation}
Here, $a=1$ for $S_2^{\prime}$ and $a=-1$ for $S_3^{\prime}$, and $p(z)=\big[12 \i \big(-9z+\sqrt{81z^2+12}\big)\big]^{1/3}$.

We note that $A_m^{(0)}$ is singular whenever $(S_m^{\prime})^2=1/3$, for which substitution into \eqref{eq:pearceyexpsol} yields two singular points, given by
\begin{equation}\label{eq:actuallocations}
z_{+}=\frac{ 2\sqrt{3}\i}{9} \qquad \text{and} \qquad z_{-}=- \frac{ 2\sqrt{3}\i}{9}.
\end{equation}
The singularities in $A_m^{(0)}$ lie at $z_{+}$ for $m=2$ and $m=3$, and at $z_{-}$ for $m=1$ and $m=2$. All of these singularities are of order $1/4$. As an example, we have that $A_1^{(0)}(z) \sim 6  C_1^{(0)}/ (-36\sqrt{3}\i[z-z_{-}])^{1/4}$.

\subsection{Integrating the singulant solutions}
In order to determine the Stokes line structure in the next section, each of \eqref{eq:valuesofS123} must be integrated to find $S_m$. We now determine these explicitly through examination of the Hamiltonian system for the Borel operator of the Pearcy equation, much like that performed by Honda \emph{et al.} \cite{honda2015virtual}.

We consider the Borel operator, $\mathscr{P}_B$, associated with the Pearcy equation \eqref{eq:pearceymain}. 
The symbol, $\sigma$, of the operator $\mathscr{P}_B$ is a complex Hamiltonian on the phase space $T^*(\mathbb{C}_w \times \mathbb{C}_z)$, where $\mathbb{C}_z$ is the physical plane and $\mathbb{C}_w$ is the Borel plane. These are given by 
\begin{equation}\label{eq:Appborelop}
  \mathscr{P}_B = \partial^3_z- \partial_w^2 \partial_z- \i z\partial_w^3 \qquad \text{and}\qquad \sigma = \xi^3 - \eta^2 \xi - \i z \eta^3,
\end{equation}
where $\xi$ and $\eta$ are local coordinates on the cotangent directions. The Hamiltonian equations arising from the natural K\"ahler structure on $T^*(\mathbb{C}_w \times \mathbb{C}_z)$,
\begin{equation}
\label{eq:hamiltonianequations}
 \dot{z}(\tau) = \partial_\xi \sigma, 
 \quad \dot{w}(\tau) = \partial_{\eta} \sigma, 
 \quad \dot{\xi}(\tau) = -\partial_{z}\sigma, 
 \quad  \dot{\eta}(\tau) = -\partial_w \sigma, 
 \quad  \sigma (\tau) = 0,
\end{equation}
may be solved to find the solutions
\begin{subequations}\label{eq:AppPartialSols}
\begin{align}
 \label{eq:AppPartialSols1}  {\eta}(\tau) &= \eta_0,\\
  \label{eq:AppPartialSols2}   {\xi}(\tau) &= \i \eta_0^3 \tau+\xi_0,\\
  \label{eq:AppPartialSols3}  {z}(\tau) &=- \eta_0^6 \tau^3+ 3 \i \xi_0 \eta_0^3 \tau^2+(3\xi_0^2-\eta_0^2)\tau +z_0,\\
  \label{eq:AppPartialSols4}  {w}(\tau) &=\frac{3 \i \eta_0^8 }{4} \tau^4+3 \eta_0^5 \xi_0 \tau^3+ \frac{\i \eta_0^4-9 \i \eta_0^2 \xi_0^2}{2}\tau^2 -(2 \eta_0\xi_0 + 3 \i \eta_0^2 z_0) \tau + w_0.
\end{align}
\end{subequations}
Here, $\eta_0$, $\xi_0$, $z_0$, and $w_0$ are constants of integration. The desired solution, $S_m$, equals $w(\tau)$ from \eqref{eq:AppPartialSols4}, which may be expressed as a function of $z$ through knowledge of $z(\tau)$ from \eqref{eq:AppPartialSols3}. The four constants of integration in \eqref{eq:AppPartialSols} are determined through imposition of the conditions
\begin{equation}\label{eq:Appboundaryconds}
\eta(0)=1, \qquad w(0)=0, \qquad z(0)=\pm \frac{2 \sqrt{3}\i}{9} ,  \qquad \sigma(0)=0.
\end{equation}
The first three of these yield $\eta_0=1$, $w_0=0$, and $z_0=\pm 2 \sqrt{3}\i/9$. The last of conditions \eqref{eq:Appboundaryconds} gives the relationship $\xi_0^3 - \xi_0 - \i z_0=0$. Since $z_0$ is known, this equation has the four solutions of $\xi_0=\pm \sqrt{3}/3$ and $\xi_0=\pm 2 \sqrt{3}/3$. The last two of these solutions give $S_m\neq 0$ at the singular point of the leading-order outer solution $A_m^{(0)}$ from \eqref{eq:pearceyO1sola}. Thus, we consider only $\xi_0=\pm \sqrt{3}/3$ for which the plus sign may be taken without loss of generality as the resultant expression for $w$ involves $\xi_0^2$ only. In defining $\hat{\tau}= - \i ( \tau - \i \xi_0)$, we substitute $\eta_0=1$, $w_0=0$, $z_0=\pm 2 \sqrt{3}\i/9$, and $\xi_0= \sqrt{3}/3$ into equations \eqref{eq:AppPartialSols3} and \eqref{eq:AppPartialSols4} to find
\begin{equation}\label{eq:AppPartialSolsnext}
   {z}(\hat{\tau}) =\i \hat{\tau}(\hat{\tau}^2-1) \qquad \text{and} \qquad {w}(\hat{\tau}) =\frac{3 \i}{4}\left(\hat{\tau}^2- 1/3 \right)^2.
\end{equation}
We recognise the first of equations \eqref{eq:AppPartialSolsnext} as the same cubic equation that governs $S_m^{\prime}$ in \eqref{eq:pearceyexpsol}, and hence $S^{\prime}_m=\hat{\tau}$. Furthermore, we may calculate $\mathrm{d}w/\mathrm{d}z = (\mathrm{d}w/\mathrm{d}\hat{\tau})(\mathrm{d}\hat{\tau}/\mathrm{d}z)=\hat{\tau}=S_m^{\prime}(z)$, which verifies that $S_m^{\prime}=w^{\prime}$. Since we have $S^{\prime}_m=\hat{\tau}$, it is not necessary to invert $z(\hat{\tau})$ to find $\hat{\tau}(z)$. Instead, we simply write our expression for $S_m$ in terms of the previously determined function $S_m^{\prime}$, which yields
\begin{equation}\label{eq:AppcheckS5}
S_m(z) =\frac{3 \i}{4}\left(\big[S_m^{\prime}(z)\big]^2- \frac{1}{3}\right)^2.
\end{equation}
Equation \eqref{eq:AppcheckS5} is the main result of this section: an explicit solution for the singulant, $S_m$, specified in terms of the function $S^{\prime}_m$ previously determined in equation \eqref{eq:valuesofS123}. 

We note that an alternative way to derive \eqref{eq:AppcheckS5} would be assume that $A^{(0)}_0(z)=\text{const} \times S_m^{-1/4}$, for which substitution into \eqref{eq:pearceyexpeqO1} yields the second-order nonlinear ODE, $12S_m S_m^{\prime \prime}-3(S_m^{\prime})^2+1=0$, satisfied by \eqref{eq:AppcheckS5}.

\subsection{Late-term divergence}
We consider a late-term ansatz of the form
\begin{equation}\label{eq:pearceyOnansatzfull}
A^{(n)}_m(z) \sim \sum_{p=0}^{\infty}B^{(p)}_m(z) \frac{\Gamma(n+\alpha-p)}{[{\chi_{m}^{}}(z)]^{n+\alpha-p}},
\end{equation}
where $\alpha$ is a constant, $\chi_m$ is the singulant, and $B_m^{(p)}$ are amplitude functions. The HOSP may be detected through analysis of the late-late term divergence of $B_m^{(p)}$ as $p \to \infty$. Substitution of \eqref{eq:pearceyOnansatzfull} into the $O(\ep^n)$ equation \eqref{eq:pearceyexpeqOepn} yields at each order of $\Gamma(n+\alpha-p)/\chi_m^{n+\alpha-p}$ the equations
\begin{subequations}\label{eq:latelateamp}
\begin{align}
\label{eq:pearcychieq}
&{\chi}_{m}^{\prime}\Big( ({\chi}_{m}^{\prime})^2+3 S_{m}^{\prime}{\chi}_{m}^{\prime}+3(S_{m}^{\prime})^2-1\Big)=0,\\
\label{eq:latelateampO1}
&\mathcal{L}\Big[B^{(0)}_m\Big] \equiv \left(3[S_m^{\prime}+\chi_m^{\prime}]^2-1\right)B^{(0)\prime}_m
+3(S_m^{\prime}+\chi_m^{\prime})(S_m^{\prime \prime}+\chi_m^{\prime \prime})B^{(0)}_m=0,\\
\label{eq:latelateampO2}
&\mathcal{L}\Big[B^{(1)}_m\Big] =3\left(S_m^{\prime}+\chi_m^{\prime} \right)B^{(0)\prime \prime}_m+3 \left(S_m^{\prime \prime}+\chi_m^{\prime \prime}\right)B^{(0)\prime}_m
+\left(S_m^{\prime \prime \prime}+\chi_m^{\prime \prime \prime}\right)B^{(0)}_m,\\
\label{eq:latelateampOp}
&\begin{aligned}
\mathcal{L}\Big[B^{(p)}_m\Big] &=3\left(S_m^{\prime}+\chi_m^{\prime} \right)B^{(p-1)\prime \prime}_{m}+3 \left(S_m^{\prime \prime}+\chi_m^{\prime \prime}\right)B^{(p-1)\prime}_{m}\\
& \quad +\left(S_m^{\prime \prime \prime}+\chi_m^{\prime \prime \prime}\right)B^{(p-1)}_{m}-B^{(p-2)\prime \prime \prime}_{m},
\end{aligned}
\end{align}
\end{subequations}
where \eqref{eq:latelateampOp} holds for $p \geq 2$, and $\mathcal{L}$ is the differential operator defined in \eqref{eq:latelateampO1}.

We begin by solving equation \eqref{eq:pearcychieq} to find $\chi_m$. Six solutions, ${\chi}_{m}^{\prime}=-3S_{m}^{\prime}/2 \pm \sqrt{4-3(S_{m}^{\prime})^2}/2$, are found in total, but two of these are unable to match with inner solutions near the associated singular points. It may be verified that these solutions for $\chi_m^{\prime}$ yield the difference between another $S_j^{\prime}$ and $S_m^{\prime}$ (for instance $\chi_1^{\prime}=S_2^{\prime}-S_1^{\prime}$).
Integration of the expressions for $\chi_m^{\prime}$ then yields the four initial singulants
\begin{equation}
\label{eq:chisolnewes2}
\left\{ \quad 
\begin{aligned}
\chi_1(z)&=S_2(z)-S_1(z), \qquad \chi_3(z)=S_2(z)-S_3(z),\\
\chi_{2a}(z) &=S_1(z)-S_2(z),\qquad  \chi_{2b}(z) =S_3(z)-S_2(z),
\end{aligned}\right.
\end{equation}
which satisfy $\chi_1(z_{-})=0$, $\chi_3(z_{+})=0$, $\chi_{2a}(z_{-})=0$, and $\chi_{2b}(z_{+})=0$. 

\subsection{Late-late-term divergence}

We now consider equations \eqref{eq:latelateamp} for the late-term amplitude functions, $B_m^{(p)}$. The first of these may be solved to find the solution
\begin{equation}\label{eq:B0sol}
B_m^{(0)}(z)=\frac{\Lambda_m^{(0)}}{\Big(3(S_m^{\prime}+\chi_m^{\prime})^2-1\Big)^{1/2}},
\end{equation}
where $\Lambda_m^{(0)}$ is a constant of integration. For $m=2$, $B_2^{(0)}$ in \eqref{eq:B0sol} is singular at the same locations as that predicted by the factorial-over-power form with $\chi_2$. However, for $m=1$ and $m=3$, \eqref{eq:B0sol}
 is singular at additional locations beyond that where the factorial-over-power form has $\chi_m(z)=0$. For example, when $m=1$, while $\chi_1$ predicts a singularity at $z_{-}$ only, the amplitude function \eqref{eq:B0sol} is also singular at $z_{+}$. This will lead to a new divergent series as $p \to \infty$, which we now study for $m=1$ and $m=3$.

In substituting a late-late-term ansatz of the form
\begin{equation}\label{eq:Bnsol}
B_m^{(p)}(z)\sim C_{m}(z) \frac{\Gamma(p+\beta)}{[\tilde{\chi}_m(z)]^{p+\beta}},
\end{equation}
into equation \eqref{eq:latelateampOp}, we find the equation $(\tilde{\chi}_m^{\prime})^2+3(S_m^{\prime}+\chi_m^{\prime})\tilde{\chi}^{\prime}_m+3(S_m^{\prime}+\chi_m^{\prime})^2-1=0$. For each value of $m$, only one solution of this quadratic equation has the correct singular behaviour after the boundary condition $\tilde{\chi}_1(z_{+})=0$ or $\tilde{\chi}_3(z_{-})=0$ is applied. This yields our late-late-term singulant as
\begin{equation}
\label{eq:chisollatelate}
\tilde{\chi}_m(z)=\left\{
\begin{aligned}
 S_{3}(z)-S_{2}(z)  \qquad \text{if $m=1$}, \\
  S_{1}(z)-S_{2}(z)  \qquad \text{if $m=3$}.
\end{aligned}\right.
\end{equation}
From the method developed in \S\ref{sec:smoothinghosp}, these will induce the HOSP, in which a new factorial-over-power component of $A_m^{(n)}$ with a singulant of $\chi_m+\tilde{\chi}_m$ will switch on across a HOSL given by $\text{Im}[\tilde{\chi}_m/\chi_m]=0$ and $\text{Re}[\tilde{\chi}_m/\chi_m] \geq 0$. The dominant components of the late-terms of the initial amplitude functions, $A_m^{(n)}$, were given in \eqref{eq:PearceyTransseirs}, and the resultant Stokes lines structure was shown in figure~\ref{fig:Pearcey}.


\providecommand{\noopsort}[1]{}
\begin{thebibliography}{10}

\bibitem{aniceto2019primer}
{\sc I.~Aniceto, G.~Ba{\c{s}}ar, and R.~Schiappa}, {\em A primer on resurgent
  transseries and their asymptotics}, Physics Reports, 809 (2019), pp.~1--135.

\bibitem{aoki1994new}
{\sc T.~Aoki}, {\em New turning points in the exact {WKB} analysis for
  higher-order ordinary differential equations}, Analyse alg{\'e}brique des
  perturbations singuli{\`e}res. I,  (1994).

\bibitem{aoki1998exact}
{\sc T.~Aoki, T.~Kawai, and Y.~Takei}, {\em On the exact {WKB} analysis for the
  third order ordinary differential equations with a large parameter}, Asian J.
  Math., 2 (1998), pp.~625--640.

\bibitem{aoki2001exact}
\leavevmode\vrule height 2pt depth -1.6pt width 23pt, {\em On the exact
  steepest descent method: {A} new method for the description of {S}tokes
  curves}, J. Math. Phys., 42 (2001), pp.~3691--3713.

\bibitem{berk1982new}
{\sc H.~L. Berk, W.~M. Nevins, and K.~V. Roberts}, {\em New {S}tokes’ line in
  {WKB} theory}, J. Math. Phys., 23 (1982), pp.~988--1002.

\bibitem{berry_1989}
{\sc M.~V. Berry}, {\em Uniform asymptotic smoothing of {S}tokes
  discontinuities}, Proc. R. Soc. Lond. A, 422 (1989), pp.~7--21.

\bibitem{berry1991hyperasymptotics}
{\sc M.~V. Berry and C.~J. Howls}, {\em Hyperasymptotics for integrals with
  saddles}, Proc. R. Soc. Lond. A, 434 (1991), pp.~657--675.

\bibitem{body2005exponential}
{\sc G.~L. Body, J.~R. King, and R.~H. Tew}, {\em Exponential asymptotics of a
  fifth-order partial differential equation}, Eur. J. Appl. Math., 16 (2005),
  pp.~647--681.

\bibitem{chapman_2005}
{\sc S.~J. Chapman and D.~B. Mortimer}, {\em Exponential asymptotics and
  {S}tokes lines in a partial differential equation}, Proc. R. Soc. Lond. A,
  461 (2005), pp.~2385--2421.

\bibitem{costin2022uniformization}
{\sc O.~Costin and G.~V. Dunne}, {\em Uniformization and constructive analytic
  continuation of {T}aylor series}, Commun. Math. Phys., 392 (2022),
  pp.~863--906.

\bibitem{crew2024resurgent}
{\sc S.~Crew and P.~H. Trinh}, {\em Resurgent aspects of applied exponential
  asymptotics}, Stud. Appl. Math., 152 (2024), pp.~974--1025.

\bibitem{dingle_book}
{\sc R.~B. Dingle}, {\em Asymptotic Expansions: Their Derivation and
  Interpretation}, Academic Press, London, 1973.

\bibitem{dorigoni2019introduction}
{\sc D.~Dorigoni}, {\em An introduction to resurgence, trans-series and alien
  calculus}, Ann. Physics, 409 (2019), p.~167914.

\bibitem{griffiths2008limiting}
{\sc S.~D. Griffiths}, {\em The limiting form of inertial instability in
  geophysical flows}, J. Fluid Mech., 605 (2008), pp.~115--143.

\bibitem{honda2007stokes}
{\sc N.~Honda}, {\em On the {S}tokes geometry of the {N}oumi-{Y}amada system},
  Algebraic, Analytic and Geometric Aspects of Complex Differential Equations
  and their Deformations. Painleve Hierarchies, 2 (2007), pp.~45--72.

\bibitem{honda2015virtual}
{\sc N.~Honda, T.~Kawai, and Y.~Takei}, {\em Virtual turning points}, vol.~4,
  Springer, 2015.

\bibitem{howls_2004}
{\sc C.~J. Howls, P.~J. Langman, and A.~B.~Olde Daalhuis}, {\em On the
  higher-order {S}tokes phenomenon}, Proc. R. Soc. Lond. A, 460 (2004),
  pp.~2285--2303.

\bibitem{iniKing}
{\sc J.~King}, {\em Exponential asymptotics, the {S}tokes phenomenon and the
  higher-order {S}tokes phenomenon in some linear partial differential
  equations}.
\newblock Issac Newton Institute Applicable resurgent asymptotics: summary
  workshop presentation 16th Dec., 2022.

\bibitem{king2001asymptotics}
{\sc J.~R. King and S.~J. Chapman}, {\em Asymptotics beyond all orders and
  {S}tokes lines in nonlinear differential-difference equations}, Eur. J. Appl.
  Math., 12 (2001), pp.~433--463.

\bibitem{lustri2019three}
{\sc C.~J. Lustri, R.~Pethiyagoda, and S.~J. Chapman}, {\em Three-dimensional
  capillary waves due to a submerged source with small surface tension}, J.
  Fluid Mech., 863 (2019), pp.~670--701.

\bibitem{mitschi2016divergent}
{\sc C.~Mitschi, D.~Sauzin, M.~Loday-Richaud, and {\'E}.~Delabaere}, {\em
  Divergent series, summability and resurgence}, Springer, 2016.

\bibitem{nemes2022dingle}
{\sc G.~Nemes}, {\em Dingle’s final main rule, {B}erry’s transition, and
  {H}owls’ conjecture}, J. Phys. A-Math. Theor., 55 (2022), p.~494001.

\bibitem{daalhuis1996hyperterminants}
{\sc A.~B. Olde~Daalhuis}, {\em Hyperterminants {I}}, J. Comp. Appl. Math., 76
  (1996), pp.~255--264.

\bibitem{daalhuis1998hyperterminants}
\leavevmode\vrule height 2pt depth -1.6pt width 23pt, {\em Hyperterminants
  {II}}, J. Comp. Appl. Math., 89 (1998), pp.~87--95.

\bibitem{daalhuis2004higher}
\leavevmode\vrule height 2pt depth -1.6pt width 23pt, {\em On higher-order
  {S}tokes phenomena of an inhomogeneous linear ordinary differential
  equation}, J. Comp. Appl. Math., 169 (2004), pp.~235--246.

\bibitem{daalhuis2009hyperasymptotics}
\leavevmode\vrule height 2pt depth -1.6pt width 23pt, {\em Hyperasymptotics and
  hyperterminants: exceptional cases}, J. Comp. Appl. Math., 233 (2009),
  pp.~555--563.

\bibitem{smoothingAdri}
\leavevmode\vrule height 2pt depth -1.6pt width 23pt, {\em Smoothing for the
  higher-order {S}tokes phenomenon}.
\newblock Issac Newton Institute Applicable resurgent asymptotics: summary
  workshop presentation 16th Dec., 2022.

\bibitem{shelton2022Hermite}
{\sc J.~Shelton, S.~J. Chapman, and P.~Trinh}, {\em Exponential asymptotics for
  a model problem of an equatorially trapped {R}ossby wave}, SIAM J. Appl.
  Math., 84 (2024), pp.~1482--1503.

\bibitem{shelton2023kelvin}
{\sc J.~Shelton, S.~Griffiths, S.~J. Chapman, and P.~H. Trinh}, {\em On the
  exponentially-small instability of the equatorial {K}elvin wave}, In
  preparation,  (2025).

\bibitem{shudo2007role}
{\sc A.~Shudo}, {\em A role of virtual turning points and new {S}tokes curves
  in {S}tokes geometry of the quantum H{\'e}non map. Algebraic analysis of
  differential equations}, Springer, 2008.

\bibitem{takei2017wkb}
{\sc Y.~Takei}, {\em {WKB} analysis and {S}tokes geometry of differential
  equations}, in Analytic, algebraic and geometric aspects of differential
  equations, Springer, 2017, pp.~263--304.

\bibitem{trinh2013new}
{\sc P.~H. Trinh and S.~J. Chapman}, {\em New gravity--capillary waves at low
  speeds. {P}art 2. {N}onlinear geometries}, J. Fluid Mech., 724 (2013),
  pp.~392--424.

\end{thebibliography}
\end{document}